\documentclass[review]{elsarticle}

\usepackage{lineno,hyperref}
\usepackage{graphicx}
\usepackage{amsmath}
\usepackage{amssymb}
\usepackage{subfigure}
\usepackage{threeparttable}
\usepackage{multirow}
\usepackage{tabularx}
\usepackage{booktabs}
\usepackage{enumerate}
\modulolinenumbers[5]

\journal{Journal of \LaTeX\ Templates}

\begin{document}

\begin{frontmatter}

\title{A minimal closed-form solution to the conic based on self-polar triangle}

\author{Yang Guo\corref{mycorrespondingauthor}}
\address{Department of Mathematics, Northeastern University, Shenyang, 110819, P.R.China}


\cortext[mycorrespondingauthor]{Corresponding author}
\ead{guoyang@mail.neu.edu.cn}


\begin{abstract}
In this paper, we use the properties of the self-polar triangle to not only show a novel method for a basic point-line enumerative problem of conics, but also present a series of closed-form solutions to the conics from all minimal configurations of points and lines in general position. These closed-form formulae may allow us to derive easily the algebraic and geometric conditions which characterize when the obtained conic is real and non-degenerate, so we propose a criterion for a non-degenerate real conic from each of all minimal configurations. The correctness of our results is validated by some examples.
\end{abstract}


\end{frontmatter}

\linenumbers

\section{Introduction}

To find the number of conics through $p$ points and tangent to $l$ lines if $p+l=5$ is a basic point-line enumerative problem of conics [1], and it can be also called as the counting problem of conics under the point-line minimal configurations. This problem has been solved by shifting the counting conics to counting the number of points in an intersection of certain subsets of the moduli space $\mathrm P\mathbb{R}^5$ [2].

Besides the enumerative problem, sometimes we care more how to compute these conics from given the point-line minimal configurations. H. D\"{o}rrie described a series of methods how to draw a conic from given the point-line minimal configurations by using Pascal's theorem (See No. 64 in the reference [3]). But these methods drew just the points on the conic satisfying the point-line minimal configurations, and could not give the analytic expression of the conic equation. Regarding for solving of the conic equation, the conventional methods were on a case-by-case basis. For example, the unique conic passing through given five points (or lines) in general position can be obtained by the vanishing of the determinant of some $6\times6$ matrix [4], or can be computed by determining the null vector of some $5\times6$ matrix [5]. Because we would get a constraint of degree 2 about the conic from a given tangent line, when a point-line mixed minimal configuration is given, the problem of solving of the conic equation can be converted to finding the common zeros of a collection of homogeneous polynomials [2]. To the best of our knowledge, on the problem of deriving the closed-form solution of the conic equation based on the self-polar triangle from all sorts of the point-line minimal correspondences, the literature is sparse. In addition, we also want especially to know that how to determine a conic satisfying a concrete given point-line minimal configuration being real and non-degenerate or not from relative position between the geometric primitives without the detailed computation [6,7].

In this paper, we use the properties of the self-polar triangle to not only show a novel method for a basic point-line enumerative problem of conics, but also propose a series of closed-form solutions to the conics satisfying all minimal configurations of points and lines in general position. Using these closed-form solutions, we propose a criterion for a non-degenerate real conic from each of all minimal configurations. Obviously, our results (whether the closed-form solutions or the criteria) may be applied in the RANSAC-based robust algorithms for estimating a conic from noisy data containing outliers [8].

This paper is organized as follows. In next section, we introduce the notations and preliminaries about the self-polar triangle. We propose a series of closed-form solutions to the conics satisfying all minimal configurations of points and lines in general position in Section 3. In Section 4, we summarize the proposed algorithms. Section 5 contains some experiments validating our theoretical results. Finally, in Section 6, conclusions are presented.

\section{Notations and Preliminaries}

In this paper, we let $\mathrm P\mathbb{R}^2$ ($\mathrm P\mathbb{C}^2$) denote the real (complex) projective plane and let $\hat {\textbf x}$ and $\textbf l$ represent respectively a normalized finite point and a line in $\mathrm P\mathbb{R}^2$, where the hat symbol $\wedge$ denote that the third component of $\hat {\textbf x}$ is one.

We say four points in $\mathrm P\mathbb{R}^2$ are in "general position", which means that no three points are collinear. It is well known that given four finite points $\hat {\textbf x}_1$, $\hat {\textbf x}_2$, $\hat {\textbf x}_3$, $\hat {\textbf x}_4$ in general position, then according to algebraic projective geometry there exists a diagonal triangle of the quadrangle $\hat {\textbf x}_1\hat {\textbf x}_2\hat {\textbf x}_3\hat {\textbf x}_4$ [9], as shown in Fig.1. We always let $\boldsymbol {\xi}_{1}$, $\boldsymbol {\xi}_{2}$, $\boldsymbol {\xi}_{3}$ denote three vertices of the diagonal triangle, and furthermore, we also stipulate that $\boldsymbol {\xi}_{1}$, $\boldsymbol {\xi}_{2}$, $\boldsymbol {\xi}_{3}$ are always computed by the following formulas:
\begin{equation}
\begin{aligned}
&\boldsymbol {\xi}_{1}=(\hat {\textbf x}_1\times\hat {\textbf x}_2)\times(\hat {\textbf x}_3\times\hat {\textbf x}_4)=det([\hat {\textbf x}_{124}])\hat {\textbf x}_{3}-det([\hat {\textbf x}_{123}])\hat {\textbf x}_{4},\\
&\boldsymbol {\xi}_{2}=(\hat {\textbf x}_1\times\hat {\textbf x}_3)\times(\hat {\textbf x}_2\times\hat {\textbf x}_4)=det([\hat {\textbf x}_{134}])\hat {\textbf x}_{2}+det([\hat {\textbf x}_{123}])\hat {\textbf x}_{4},\\
&\boldsymbol {\xi}_{3}=(\hat {\textbf x}_1\times\hat {\textbf x}_4)\times(\hat {\textbf x}_2\times\hat {\textbf x}_3)=det([\hat {\textbf x}_{124}])\hat {\textbf x}_{3}-det([\hat {\textbf x}_{134}])\hat {\textbf x}_{2},
\end{aligned}
\end{equation}
where $\times$ denotes the cross product and $[\hat {\textbf x}_{ijk}]$ expresses the $3\times3$ matrix $[\hat {\textbf x}_i, \hat {\textbf x}_j, \hat {\textbf x}_k]$. $det(*)$ denotes the determinant function.

We use $[\hat {\textbf x}_{ijk}]^T$ and $[\hat {\textbf x}_{ijk}]^{-1}$ to denote the transpose matrix and the inverse matrix of $[\hat {\textbf x}_{ijk}]$ respectively. The symbol $\simeq$ denotes equality up to scale.

We assume that the reader is very familiar with the self-polar triangle for a non-degenerate conic. However, for the ease of the presentation, we nevertheless list an important result without elaboration.

Given four points $\hat {\textbf x}_1$, $\hat {\textbf x}_2$, $\hat {\textbf x}_3$, $\hat {\textbf x}_4\in\mathrm P\mathbb{R}^2$ in general position. Then the diagonal triangle $\boldsymbol {\xi}_{1}\boldsymbol {\xi}_{2}\boldsymbol {\xi}_{3}$ of the quadrangle $\hat {\textbf x}_1\hat {\textbf x}_2\hat {\textbf x}_3\hat {\textbf x}_4$ is self-polar such that
\begin{equation*}
\begin{aligned}
\boldsymbol {\xi}_{i}^T\textbf C\boldsymbol {\xi}_{j}=0, i\neq j\in\{1, 2, 3\},
\end{aligned}
\end{equation*}
for any non-degenerate conic $\textbf C$ passing through each of $\hat {\textbf x}_1$, $\hat {\textbf x}_2$, $\hat {\textbf x}_3$, $\hat {\textbf x}_4$.

\section{A minimal closed-form solution to the conic}

In this section, first of all, we derive an analytic solution with one degree of freedom to the non-degenerate conic based on self-polar triangle from four-point configuration. Further, we present a series of closed-form solutions to the conics from all minimal configurations of points and lines in general position. In order to do that, we will need the following lemmas.\\

\begin{figure}
\begin{center}
  \includegraphics[width=0.6\textwidth]{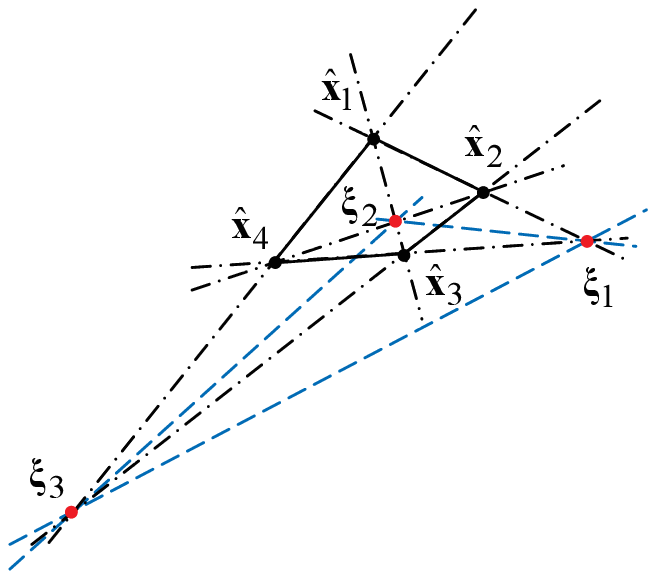}
\caption{$\boldsymbol {\xi}_1$, $\boldsymbol {\xi}_2$, $\boldsymbol {\xi}_3$ are three intersection points of three pairs of opposite sides of the quadrangle $\hat {\textbf x}_1\hat {\textbf x}_2\hat {\textbf x}_3\hat {\textbf x}_4$. The triangle $\boldsymbol {\xi}_1\boldsymbol {\xi}_2\boldsymbol {\xi}_3$ is called the diagonal triangle of the quadrangle $\hat {\textbf x}_1\hat {\textbf x}_2\hat {\textbf x}_3\hat {\textbf x}_4$.}
\label{fig:1}       
\end{center}
\end{figure}

\noindent\textbf {Lemma 1.} \emph{Given a triangle whose three vertices are} $\boldsymbol {\xi}_{1}, \boldsymbol {\xi}_{2}, \boldsymbol {\xi}_{3}\in\!\mathrm P\mathbb{R}^2$. \emph{Then the set}\\
\begin{equation*}
\begin{aligned}
S_{\boldsymbol {\xi}_{123}}=\{\ \textbf C\ |\ \boldsymbol {\xi}_{i}^T\textbf C\boldsymbol {\xi}_{j}=0, i\neq j\in\{1, 2, 3\}, \textbf C\in\mathbb{R}^{3\times3}\},
\end{aligned}
\end{equation*}
\emph{is a three-dimensional vector space over} $\mathbb{R}$.\\

\noindent\textbf {Proof.} The proof that the set $S_{\boldsymbol {\xi}_{123}}$ is a vector space over $\mathbb{R}$ is very straightforward according to the definition of vector space [10]. In the follows, we will prove that the dimension of $S_{\boldsymbol {\xi}_{123}}$ is three.

Since $\boldsymbol {\xi}_{j}^T\textbf C\boldsymbol {\xi}_{i}=0$, $\boldsymbol {\xi}_{k}^T\textbf C\boldsymbol {\xi}_{i}=0$, $i\neq j\neq k\in\{1, 2, 3\}$, we have $\textbf C\boldsymbol {\xi}_{i}=\lambda_i(\boldsymbol {\xi}_{j}\times\boldsymbol {\xi}_{k})$, where $\lambda_i$ is an unknown real number. Further, we may get
\begin{equation*}
\begin{aligned}
\textbf C[\boldsymbol {\xi}_{123}]=[\boldsymbol {\xi}_{2}\times\boldsymbol {\xi}_{3}, \boldsymbol {\xi}_{1}\times\boldsymbol {\xi}_{3}, \boldsymbol {\xi}_{1}\times\boldsymbol {\xi}_{2}]Diag(\lambda_1, \lambda_2, \lambda_3).
\end{aligned}
\end{equation*}
Since
\begin{equation*}
\begin{aligned}
\ [\boldsymbol {\xi}_{123}]^{-T}=\frac{1}{det([\boldsymbol {\xi}_{123}])}[\boldsymbol {\xi}_{2}\!\times\!\boldsymbol {\xi}_{3}, -\boldsymbol {\xi}_{1}\!\times\!\boldsymbol {\xi}_{3}, \boldsymbol {\xi}_{1}\!\times\!\boldsymbol {\xi}_{2}],
\end{aligned}
\end{equation*}
we may derive that
\begin{equation*}
\begin{aligned}
\textbf C=det([\boldsymbol {\xi}_{123}])[\boldsymbol {\xi}_{123}]^{-T}Diag(\lambda_1, -\lambda_2, \lambda_3)[\boldsymbol {\xi}_{123}]^{-1},
\end{aligned}
\end{equation*}
where $\lambda_i$, $i=1, 2, 3$ are unknown real numbers.

Now, let
\begin{equation*}
\begin{aligned}
&\textbf B_1=det([\boldsymbol {\xi}_{123}])[\boldsymbol {\xi}_{123}]^{-T}Diag(1, 0, 0)[\boldsymbol {\xi}_{123}]^{-1},\\
&\textbf B_2=det([\boldsymbol {\xi}_{123}])[\boldsymbol {\xi}_{123}]^{-T}Diag(0, 1, 0)[\boldsymbol {\xi}_{123}]^{-1},\\
&\textbf B_3=det([\boldsymbol {\xi}_{123}])[\boldsymbol {\xi}_{123}]^{-T}Diag(0, 0, 1)[\boldsymbol {\xi}_{123}]^{-1}.
\end{aligned}
\end{equation*}
Obviously, each vector in the vector space $S_{\boldsymbol {\xi}_{123}}$ is a linear combination of $\textbf B_1$, $\textbf B_2$, $\textbf B_3$ and furthermore, $\textbf B_1$, $\textbf B_2$, $\textbf B_3$ are linearly independent, so $\{\textbf B_1$, $\textbf B_2$, $\textbf B_3\}$ is a basis of the vector space $S_{\boldsymbol {\xi}_{123}}$. That means that the vector space $S_{\boldsymbol {\xi}_{123}}$ is three-dimensional. This completes the proof of the lemma. \hfill {$\Box$}\\

From the proof of Lemma 1, we know that each vector in the vector space $S_{\boldsymbol {\xi}_{123}}$ is a real symmetric matrix.\\

\noindent\textbf {Lemma 2.} \emph{Given four points} $\hat {\textbf x}_{1}, \hat {\textbf x}_{2}, \hat {\textbf x}_{3}, \hat {\textbf x}_{4}\in\!\mathrm P\mathbb{R}^2$ \emph{in general position}. \emph{Let} $\boldsymbol {\xi}_{1}, \boldsymbol {\xi}_{2}, \boldsymbol {\xi}_{3}$ \emph{are three vertices of the diagonal triangle of the quadrangle} $\hat {\textbf x}_1\hat {\textbf x}_2\hat {\textbf x}_3\hat {\textbf x}_4$. \emph{Now given one more finite point} $\hat {\textbf x}_5\in\!\mathrm P\mathbb{R}^2$, \emph{then} $\hat {\textbf x}_5$ \emph{lies on some side of the quadrangle} $\hat {\textbf x}_1\hat {\textbf x}_2\hat {\textbf x}_3\hat {\textbf x}_4$ \emph{if and only if the absolute values of at least two components are equal in the vector} $[\boldsymbol {\xi}_{123}]^{-1}\hat {\textbf x}_5$.\\

\noindent\textbf {Proof.} First, according to (1), we may derive
\begin{equation*}
\begin{aligned}
&\ [\boldsymbol {\xi}_{123}]^{-1}\hat {\textbf x}_1\simeq(-1,1,-1)^T,\quad [\boldsymbol {\xi}_{123}]^{-1}\hat {\textbf x}_2\simeq(1,1,-1)^T,\\
&\ [\boldsymbol {\xi}_{123}]^{-1}\hat {\textbf x}_3\simeq(1,1,1)^T,\qquad [\boldsymbol {\xi}_{123}]^{-1}\hat {\textbf x}_4\simeq(-1,1,1)^T.
\end{aligned}
\end{equation*}
Further, we may know that if the point $\hat {\textbf x}_5$ lies on the line through $\hat {\textbf x}_i, \hat {\textbf x}_j$, $i\neq j\in\{1, 2, 3, 4\}$, namely $\hat {\textbf x}_5=k_i\hat {\textbf x}_i+k_j\hat {\textbf x}_j$, where $k_i, k_j\in\mathbb{R}$ satisfying $k_i+k_j=1$, then the absolute values of at least two components must be equal in the vector $[\boldsymbol {\xi}_{123}]^{-1}\hat {\textbf x}_5$.

Conversely, if the absolute values of at least two components are equal in the vector $[\boldsymbol {\xi}_{123}]^{-1}\hat {\textbf x}_5$, then $[\boldsymbol {\xi}_{123}]^{-1}\hat {\textbf x}_5$ must be a linear combination of some two vectors among $[\boldsymbol {\xi}_{123}]^{-1}\hat {\textbf x}_i$, $i=1,2,3,4$. That means that the point $\hat {\textbf x}_5$ must lie on some side of the quadrangle $\hat {\textbf x}_1\hat {\textbf x}_2\hat {\textbf x}_3\hat {\textbf x}_4$. This completes the proof of the lemma. \hfill {$\Box$}\\

\subsection{Four Points In General Position}

\noindent\textbf {Theorem 1.} \emph{Given four points} $\hat {\textbf x}_{1}, \hat {\textbf x}_{2}, \hat {\textbf x}_{3}, \hat {\textbf x}_{4}\in\mathrm P\mathbb{R}^2$ \emph{in general position}. \emph{Let} $\boldsymbol {\xi}_{1}, \boldsymbol {\xi}_{2}, \boldsymbol {\xi}_{3}$ \emph{are three vertices of the diagonal triangle of the quadrangle} $\hat {\textbf x}_1\hat {\textbf x}_2\hat {\textbf x}_3\hat {\textbf x}_4$. \emph{Then}
\begin{equation}
\begin{aligned}
\textbf C=[\boldsymbol {\xi}_{123}]^{-T}Diag(-s, s-1, 1)[\boldsymbol {\xi}_{123}]^{-1},
\end{aligned}
\end{equation}
\emph{is a non-degenerate conic with one degree of freedom such that} $\hat {\textbf x}_i^T\textbf C\hat {\textbf x}_i=0$, $i=1,2,3,4$, \emph{where} $s\in\mathbb{R}\backslash \{0, 1\}$.\\

\noindent\textbf {Proof.} From algebraic projective geometry, we know that $\boldsymbol {\xi}_{1}, \boldsymbol {\xi}_{2}, \boldsymbol {\xi}_{3}$ form a self-polar triangle such that
\begin{equation*}
\begin{aligned}
\boldsymbol {\xi}_{i}^T\textbf C\boldsymbol {\xi}_{j}=0, i\neq j\in\{1, 2, 3\},
\end{aligned}
\end{equation*}
for any non-degenerate conic $\textbf C$ passing through $\hat {\textbf x}_{1}$, $\hat {\textbf x}_{2}$, $\hat {\textbf x}_{3}$, $\hat {\textbf x}_{4}$. So according to the proof of Lemma 1, we have
\begin{equation*}
\begin{aligned}
\textbf C\simeq[\boldsymbol {\xi}_{123}]^{-T}Diag(\lambda_1, -\lambda_2, \lambda_3)[\boldsymbol {\xi}_{123}]^{-1},
\end{aligned}
\end{equation*}
where $\lambda_1$, $\lambda_2$, $\lambda_3$ are unknown non-zero real numbers.

Further, according to the proof of Lemma 2 and $\hat {\textbf x}_i^T\textbf C\hat {\textbf x}_i=0$, $i=1,2,3,4$, we have $\lambda_1-\lambda_2+\lambda_3=0$. Since $\textbf C$ can be determined up to scale, we assert that a non-degenerate conic passing through four points in general position may be determined with one degree of freedom. This completes the proof of the theorem. \hfill {$\Box$}\\

\subsection{Five Points In General Position}

In this paper, we say five points are in "general position", which means that no three points are collinear.\\

\noindent\textbf {Theorem 2.} \emph{Given five points} $\hat {\textbf x}_{1}, \hat {\textbf x}_{2}, \hat {\textbf x}_{3}, \hat {\textbf x}_{4}, \hat {\textbf x}_{5}\in\mathrm P\mathbb{R}^2$ \emph{in general position}. \emph{Let} $\boldsymbol {\xi}_{1}, \boldsymbol {\xi}_{2}, \boldsymbol {\xi}_{3}$ \emph{are three vertices of the diagonal triangle of the quadrangle} $\hat {\textbf x}_1\hat {\textbf x}_2\hat {\textbf x}_3\hat {\textbf x}_4$. \emph{Then}
\begin{equation}
\begin{aligned}
&\textbf C=[\boldsymbol {\xi}_{123}]^{-T}Diag(\beta_{\textbf x_53}^2-\beta_{\textbf x_52}^2, \beta_{\textbf x_51}^2-\beta_{\textbf x_53}^2, \beta_{\textbf x_52}^2-\beta_{\textbf x_51}^2)[\boldsymbol {\xi}_{123}]^{-1},
\end{aligned}
\end{equation}
\emph{is a unique non-degenerate conic such that} $\hat {\textbf x}_i^T\textbf C\hat {\textbf x}_i=0$, $i=1,2,3,4,5$, \emph{where} $(\beta_{\textbf x_51}, \beta_{\textbf x_52}, \beta_{\textbf x_53})^T=[\boldsymbol {\xi}_{123}]^{-1}\hat {\textbf x}_5$.\\

\noindent\textbf {Proof.} Using $\hat {\textbf x}_{1}, \hat {\textbf x}_{2}, \hat {\textbf x}_{3}, \hat {\textbf x}_{4}$ and Theorem 1, we may get a non-degenerate conic with one degree of freedom. Now, let $(\beta_{\textbf x_51}, \beta_{\textbf x_52}, \beta_{\textbf x_53})^T=[\boldsymbol {\xi}_{123}]^{-1}\hat {\textbf x}_5$, according to Lemma 2, we know that the absolute values of any two among $\beta_{\textbf x_51}$, $\beta_{\textbf x_52}$, $\beta_{\textbf x_53}$ must be not equal. So according to $\hat {\textbf x}_5^T\textbf C\hat {\textbf x}_5=0$ and (2), we may obtain
\begin{equation*}
\begin{aligned}
s=\frac{\beta_{\textbf x_52}^2-\beta_{\textbf x_53}^2}{\beta_{\textbf x_52}^2-\beta_{\textbf x_51}^2}.
\end{aligned}
\end{equation*}
Further, we substitute it in (2), then we may get (3). This completes the proof of the theorem. \hfill {$\Box$}\\

\subsection{Four-point + One-line In General Position}
In this paper, we say four points $\hat {\textbf x}_1$, $\hat {\textbf x}_2$, $\hat {\textbf x}_3$, $\hat {\textbf x}_4$ and one line $\textbf l$ are in "general position", which means that $\hat {\textbf x}_1$, $\hat {\textbf x}_2$, $\hat {\textbf x}_3$, $\hat {\textbf x}_4$ are in general position and $\textbf l$ is neither a side or a diagonal line of the quadrangle $\hat {\textbf x}_1\hat {\textbf x}_2\hat {\textbf x}_3\hat {\textbf x}_4$ nor a side of the diagonal triangle of the quadrangle $\hat {\textbf x}_1\hat {\textbf x}_2\hat {\textbf x}_3\hat {\textbf x}_4$.\\

\noindent\textbf {Lemma 3.} \emph{Given four points} $\hat {\textbf x}_{1}, \hat {\textbf x}_{2}, \hat {\textbf x}_{3}, \hat {\textbf x}_{4}\in\mathrm P\mathbb{R}^2$ \emph{and a line} $\textbf l\in\mathrm P\mathbb{R}^2$ \emph{in general position}. \emph{Let} $\boldsymbol {\xi}_{1}, \boldsymbol {\xi}_{2}, \boldsymbol {\xi}_{3}$ \emph{are three vertices of the diagonal triangle of the quadrangle} $\hat {\textbf x}_1\hat {\textbf x}_2\hat {\textbf x}_3\hat {\textbf x}_4$. \emph{If} $\boldsymbol {\xi}_{i}^T\textbf l=0$, \emph{then} $(\boldsymbol {\xi}_{j}^T\textbf l)^2\neq (\boldsymbol {\xi}_{k}^T\textbf l)^2$, \emph{where} $i\neq j\neq k\in\{1,2,3\}$.\\

\noindent\textbf {Proof.} We prove $i=1$, $j=2$, $k=3$, and the proofs of other cases are similar.

First, according to (1) and $\boldsymbol {\xi}_{1}^T\textbf l=0$, we have
\begin{equation}
\begin{aligned}
det([\hat {\textbf x}_{124}])\hat {\textbf x}_3^T\textbf l-det([\hat {\textbf x}_{123}])\hat {\textbf x}_4^T\textbf l=-det([\hat {\textbf x}_{234}])\hat {\textbf x}_1^T\textbf l+det([\hat {\textbf x}_{134}])\hat {\textbf x}_2^T\textbf l=0.
\end{aligned}
\end{equation}

Second, if
\begin{equation*}
\begin{aligned}
(\boldsymbol {\xi}_{2}^T\textbf l)^2-(\boldsymbol {\xi}_{3}^T\textbf l)^2=(\boldsymbol {\xi}_{2}^T\textbf l+\boldsymbol {\xi}_{3}^T\textbf l)(\boldsymbol {\xi}_{2}^T\textbf l-\boldsymbol {\xi}_{3}^T\textbf l)=0,
\end{aligned}
\end{equation*}
then
\begin{equation}
\begin{aligned}
\boldsymbol {\xi}_{2}^T\textbf l+\boldsymbol {\xi}_{3}^T\textbf l=det([\hat {\textbf x}_{124}])\hat {\textbf x}_3^T\textbf l+det([\hat {\textbf x}_{123}])\hat {\textbf x}_4^T\textbf l=0,
\end{aligned}
\end{equation}
or
\begin{equation}
\begin{aligned}
\boldsymbol {\xi}_{2}^T\textbf l-\boldsymbol {\xi}_{3}^T\textbf l=2det([\hat {\textbf x}_{134}])\hat {\textbf x}_2^T\textbf l=0.
\end{aligned}
\end{equation}

According to (4) and (5) or (4) and (6), we may derive that the line $\textbf l$ must pass through $\hat {\textbf x}_3$ and $\hat {\textbf x}_4$ or through $\hat {\textbf x}_1$ and $\hat {\textbf x}_2$. That is a contradiction with the fact that $\hat {\textbf x}_{1}, \hat {\textbf x}_{2}, \hat {\textbf x}_{3}, \hat {\textbf x}_{4}$ and $\textbf l$ are in general position.  This completes the proof of the lemma. \hfill {$\Box$}\\

\noindent\textbf {Theorem 3.} \emph{Given four points} $\hat {\textbf x}_{1}, \hat {\textbf x}_{2}, \hat {\textbf x}_{3}, \hat {\textbf x}_{4}\in\mathrm P\mathbb{R}^2$ \emph{and a line} $\textbf l\in\mathrm P\mathbb{R}^2$ \emph{in general position}. \emph{Let} $\boldsymbol {\xi}_{1}, \boldsymbol {\xi}_{2}, \boldsymbol {\xi}_{3}$ \emph{are three vertices of the diagonal triangle of the quadrangle} $\hat {\textbf x}_{1}\hat {\textbf x}_{2}\hat {\textbf x}_{3}\hat {\textbf x}_{4}$. \emph{Then the problem of determining a non-degenerate conic such that} $\hat {\textbf x}_i^T\textbf C\hat {\textbf x}_i=0$, $i=1,2,3,4$ \emph{and} $\textbf l^T\textbf C^{-1}\textbf l\!=\!0$ \emph{has at most two solutions in} $\mathrm P\mathbb{C}^2$ \emph{as follows}:
\begin{equation}
\begin{aligned}
\textbf C_i=[\boldsymbol {\xi}_{123}]^{-T}Diag(-s_i, s_i-1, 1)[\boldsymbol {\xi}_{123}]^{-1}, i=1,2,
\end{aligned}
\end{equation}
\emph{where} $s_i\neq0,1$ \emph{is a root of the equation}
\begin{equation}
\begin{aligned}
(\boldsymbol {\xi}_{3}^T\textbf l)^2s^2\!-\!((\boldsymbol {\xi}_{1}^T\textbf l)^2\!-\!(\boldsymbol {\xi}_{2}^T\textbf l)^2\!+\!(\boldsymbol {\xi}_{3}^T\textbf l)^2)s\!+\!(\boldsymbol {\xi}_{1}^T\textbf l)^2\!=\!0.\\
\end{aligned}
\end{equation}

\noindent\textbf {Proof.} Using $\hat {\textbf x}_{1}, \hat {\textbf x}_{2}, \hat {\textbf x}_{3}, \hat {\textbf x}_{4}$ and Theorem 1, we may get a non-degenerate conic with one degree of freedom. Further, according to the constraint $\textbf l^T\textbf C^{-1}\textbf l=0$, we may get (8). If (8) has two distinct roots $s_1, s_2$ satisfying $s_1, s_2\neq 0, 1$, then the problem of determining a non-degenerate conic must have two distinct solutions. This completes the proof of the theorem. \hfill {$\Box$}\\

By the detailed derivation, the discriminant of (8) may be represented as follows:
\begin{equation}
\begin{aligned}
\Delta&=((\boldsymbol {\xi}_{1}^T\textbf l)^2\!-\!(\boldsymbol {\xi}_{2}^T\textbf l)^2\!+\!(\boldsymbol {\xi}_{3}^T\textbf l)^2)^2-4(\boldsymbol {\xi}_{1}^T\textbf l)^2(\boldsymbol {\xi}_{3}^T\textbf l)^2\\
&=16det([\hat {\textbf x}_{123}])det([\hat {\textbf x}_{124}])det([\hat {\textbf x}_{134}])det([\hat {\textbf x}_{234}])(\hat {\textbf x}_1^T\textbf l)(\hat {\textbf x}_2^T\textbf l)(\hat {\textbf x}_3^T\textbf l)(\hat {\textbf x}_4^T\textbf l).
\end{aligned}
\end{equation}
On one hand, obviously, when the line $\textbf l$ passes through some one of $\hat {\textbf x}_1$, $\hat {\textbf x}_2$, $\hat {\textbf x}_3$ and $\hat {\textbf x}_4$, we have $\Delta=0$. That means that $s$ in (8) has a unique real root as follows:
\begin{equation}
\begin{aligned}
s=\frac{(\boldsymbol {\xi}_{1}^T\textbf l)^2-(\boldsymbol {\xi}_{2}^T\textbf l)^2+(\boldsymbol {\xi}_{3}^T\textbf l)^2}{2(\boldsymbol {\xi}_{3}^T\textbf l)^2}.
\end{aligned}
\end{equation}
It is worth noting that under this case, $\boldsymbol {\xi}_{i}^T\textbf l$, $i=1,2,3$ must be not equal to zero, otherwise $\textbf l$ is a side or diagonal line of the quadrangle $\hat {\textbf x}_1\hat {\textbf x}_2\hat {\textbf x}_3\hat {\textbf x}_4$. Consequently, $s$ in (8) must be not $0$ or $1$ (otherwise $\boldsymbol {\xi}_{1}^T\textbf l=0$ or $\boldsymbol {\xi}_{2}^T\textbf l=0$).

On the other hand, when the line $\textbf l$ passes through some one of $\boldsymbol {\xi}_{1}$, $\boldsymbol {\xi}_{2}$ and $\boldsymbol {\xi}_{3}$, the problem of determining a non-degenerate conic has also only a unique real solution. For example, when $\boldsymbol {\xi}_{3}^T\textbf l=0$, we have
\begin{equation}
\begin{aligned}
s=\frac{(\boldsymbol {\xi}_{1}^T\textbf l)^2}{(\boldsymbol {\xi}_{1}^T\textbf l)^2-(\boldsymbol {\xi}_{2}^T\textbf l)^2}.
\end{aligned}
\end{equation}
When $\boldsymbol {\xi}_{1}^T\textbf l=0$, we have
\begin{equation*}
\begin{aligned}
s_1=1-\frac{(\boldsymbol {\xi}_{2}^T\textbf l)^2}{(\boldsymbol {\xi}_{3}^T\textbf l)^2}, \quad\quad s_2=0.
\end{aligned}
\end{equation*}
Obviously, only $s_1$ can give a non-degenerate real conic.
When $\boldsymbol {\xi}_{2}^T\textbf l=0$, we have
\begin{equation*}
\begin{aligned}
s_1=\frac{(\boldsymbol {\xi}_{1}^T\textbf l)^2}{(\boldsymbol {\xi}_{3}^T\textbf l)^2}, \quad\quad s_2=1.
\end{aligned}
\end{equation*}
Certainly, only $s_1$ can give a non-degenerate real conic.

The following proposition summarizes our findings.\\

\noindent\textbf {Proposition 1.} \emph{Given four points} $\hat {\textbf x}_1$, $\hat {\textbf x}_2$, $\hat {\textbf x}_3$, $\hat {\textbf x}_4\in\mathrm P\mathbb{R}^2$ \emph{and a line} $\textbf l\in\mathrm P\mathbb{R}^2$ \emph{in general position. Then the problem of determining a non-degenerate conic has a unique real solution if and only if the line} $\textbf l$ \emph{passes through some one of} $\hat {\textbf x}_1$, $\hat {\textbf x}_2$, $\hat {\textbf x}_3$, $\hat {\textbf x}_4$, $\boldsymbol {\xi}_{1}$, $\boldsymbol {\xi}_{2}$, $\boldsymbol {\xi}_{3}$, \emph{where} $\boldsymbol {\xi}_{i}$, $i=1,2,3$ \emph{are three vertices of the diagonal triangle of the quadrangle} $\hat {\textbf x}_1\hat {\textbf x}_2\hat {\textbf x}_3\hat {\textbf x}_4$.\\

Except for the case described in Proposition 1, all other cases give two distinct solutions for determining a non-degenerate conic problem in $\mathrm P\mathbb{C}^2$. What we are interesting is that under what configuration the problem of determining a non-degenerate conic has two distinct real solutions. In order to explore this problem, we need the following concept.\\

\noindent\textbf {Definition 1.} Given three non-collinear finite points $\hat {\textbf p}_1$, $\hat {\textbf p}_2$, $\hat {\textbf p}_3\in\mathrm P\mathbb{R}^2$, if the movement direction of $\hat {\textbf p}_1\rightarrow\hat {\textbf p}_2\rightarrow\hat {\textbf p}_3$ is clockwise, then we call the determinant $det([\hat {\textbf p}_{123}])$ clockwise. Otherwise we call it anticlockwise.\\

It is well known that if $det([\hat {\textbf p}_{123}])$ is anticlockwise, then $det([\hat {\textbf p}_{123}])>0$, otherwise $det([\hat {\textbf p}_{123}])<0$ [11]. Further, we may infer that given four points $\hat {\textbf x}_{1}$, $\hat {\textbf x}_{2}$, $\hat {\textbf x}_{3}$, $\hat {\textbf x}_{4}$ in general position in $\mathrm P\mathbb{R}^2$, then
\begin{equation*}
\begin{aligned}
det([\hat {\textbf x}_{123}])det([\hat {\textbf x}_{124}])det([\hat {\textbf x}_{134}])det([\hat {\textbf x}_{234}])
\end{aligned}
\end{equation*}
must be invariant with respect to permutation methods of the four points.\\

In addition, we also know that given a non-infinite line $\textbf l$, all finite points not lying on $\textbf l$ may be divided into two parts by $\textbf l$ according to the sign of $\hat {\textbf x}^{T}\textbf l$, namely one side is "+" and another side is "-", where $\hat {\textbf x}$ is a normalized finite point.\\

According to the above analysis, we may deduce that under what configuration the problem of determining a non-degenerate conic has two distinct real solutions. More formally, we have the following proposition.\\

\noindent\textbf {Proposition 2.} \emph{Given four points} $\hat {\textbf x}_{1}, \hat {\textbf x}_{2}, \hat {\textbf x}_{3}, \hat {\textbf x}_{4}\in\mathrm P\mathbb{R}^2$ \emph{and a line} $\textbf l\in\mathrm P\mathbb{R}^2$ \emph{in general position. Then the problem of determining a non-degenerate conic has two distinct real solutions if and only if} $\textbf l$ \emph{does not pass through any one of $\hat {\textbf x}_{1}$, $\hat {\textbf x}_{2}$, $\hat {\textbf x}_{3}$, $\hat {\textbf x}_{4}$ and three vertices of the diagonal triangle of the quadrangle $\hat {\textbf x}_{1}\hat {\textbf x}_{2}\hat {\textbf x}_{3}\hat {\textbf x}_{4}$ and furthermore one of the following conditions holds}:\\

\par\setlength\parindent{0.5em}(i) \emph{all} $det([\hat {\textbf x}_{123}])$, $det([\hat {\textbf x}_{124}])$, $det([\hat {\textbf x}_{134}])$, $det([\hat {\textbf x}_{234}])$ \emph{have same directions or some two are clockwise and the remaining two are anticlockwise} + \emph{all four points lie on the same side of the line or some two points lie on a side and the remaining two lie on another side};\\

\par\setlength\parindent{0.2em}(ii) \emph{some three of} $det([\hat {\textbf x}_{123}])$, $det([\hat {\textbf x}_{124}])$, $det([\hat {\textbf x}_{134}])$, $det([\hat {\textbf x}_{234}])$ \emph{are same directions and remaining one is anti-direction} + \emph{some three points lie on the same side of the line and remaining one lies on another side}.\\

\subsection{Three-point + Two-line In General Position}
In this paper, we say three points $\hat {\textbf x}_1$, $\hat {\textbf x}_2$, $\hat {\textbf x}_3$ and two lines $\textbf l_1$, $\textbf l_2$ are in "general position", which means that $\hat {\textbf x}_1$, $\hat {\textbf x}_2$, $\hat {\textbf x}_3$ are noncollinear and neither $\textbf l_1$ nor $\textbf l_2$ is a side of the triangle $\hat {\textbf x}_1\hat {\textbf x}_2\hat {\textbf x}_3$ and moreover the intersection point of $\textbf l_1$ with $\textbf l_2$ is not a vertex of the triangle $\hat {\textbf x}_1\hat {\textbf x}_2\hat {\textbf x}_3$.\\

\noindent\textbf {Theorem 4.} \emph{Given three points} $\hat {\textbf x}_{1}$, $\hat {\textbf x}_{2}$, $\hat {\textbf x}_{3}\in\mathrm P\mathbb{R}^2$ \emph{and two lines} $\textbf l_1$, $\textbf l_2\in\mathrm P\mathbb{R}^2$ \emph{in general position}. \emph{If} $\hat {\textbf x}_i^T\textbf l_1\neq0$, $\hat {\textbf x}_i^T\textbf l_2\neq0$, $i=1,2,3$ \emph{and the point} $\textbf l_1\times\textbf l_2$ \emph{does not lie on any side of the triangle} $\hat {\textbf x}_1\hat {\textbf x}_2\hat {\textbf x}_3$, \emph{then the problem of determining a non-degenerate conic must have four distinct solutions}.\\

\noindent\textbf {Proof.} In this case, we let $\textbf x_4=t\hat {\textbf x}_1+\textbf p$, where $\textbf p=\textbf l_1\times\textbf l_2$. According to (1) and (8), we have
\begin{equation}
\begin{aligned}
(s,t,1)\textbf C_{s,t}^{(i)}\left[\begin{array}{c}
    s \\  
    t \\  
    1 \\  
  \end{array}
\right]=0, \quad i=1,2,\\
\end{aligned}
\end{equation}
where
\begin{equation*}
\begin{aligned}
\textbf C_{s,t}^{(i)}=\left[\begin{array}{ccc}
    c_{i1} & c_{i2} & c_{i3} \\  
    c_{i2} & c_{i4} & c_{i5} \\  
    c_{i3} & c_{i5} & c_{i6} \\  
  \end{array}
\right],\quad  i=1,2,\\
\end{aligned}
\end{equation*}
and
\begin{equation}
\begin{aligned}
&c_{i1}=det([\textbf p, \hat {\textbf x}_2, \hat {\textbf x}_3])^2(\hat {\textbf x}_{1}^T\textbf l_i)^2,\\
&c_{i2}=(det([\hat {\textbf x}_1, \hat {\textbf x}_2, \textbf p])\hat {\textbf x}_3^T\textbf l_i-det([\hat {\textbf x}_1, \textbf p, \hat {\textbf x}_3])\hat {\textbf x}_2^T\textbf l_i)det([\hat {\textbf x}_{123}])\hat {\textbf x}_1^T\textbf l_i,\\
&c_{i3}=det([\textbf p, \hat {\textbf x}_2, \hat {\textbf x}_3])(det([\hat {\textbf x}_1, \hat {\textbf x}_2, \textbf p])(\hat {\textbf x}_{1}^T\textbf l_i)(\hat {\textbf x}_3^T\textbf l_i),\\
&c_{i4}=det([\hat {\textbf x}_{123}])^2(\hat {\textbf x}_1^T\textbf l_i)^2,\\
&c_{i5}=-det([\hat {\textbf x}_{123}])det([\hat {\textbf x}_1, \hat {\textbf x}_2, \textbf p])(\hat {\textbf x}_1^T\textbf l_i)(\hat {\textbf x}_3^T\textbf l_i),\\
&c_{i6}=det([\hat {\textbf x}_1, \hat {\textbf x}_2, \textbf p])^2(\hat {\textbf x}_3^T\textbf l_i)^2.
\end{aligned}
\end{equation}

Further, we may derive that
\begin{equation}
\begin{aligned}
&det(\textbf C_{s,t}^{(i)})\!=\!-4(\hat {\textbf x}_1^T\textbf l_i)^2(\hat {\textbf x}_2^T\textbf l_i)^2(\hat {\textbf x}_3^T\textbf l_i)^2\\
&\quad\quad\quad\quad\quad\quad\quad\quad\quad det([\hat {\textbf x}_{123}])^2det([\hat {\textbf x}_1, \textbf p, \hat {\textbf x}_3])^2det([\hat {\textbf x}_1, \hat {\textbf x}_2, \textbf p])^2\neq0.\\
\end{aligned}
\end{equation}
That means that the problem is converted to solve the intersection points of two non-degenerate conics in variables $s$, $t$.

By careful computation, we may know that
\begin{equation}
\begin{aligned}
&\lambda_1=\frac{(\hat {\textbf x}_1^T\textbf l_2)(\hat {\textbf x}_2^T\textbf l_2)}{(\hat {\textbf x}_1^T\textbf l_1)(\hat {\textbf x}_2^T\textbf l_1)}, \quad \quad \lambda_2=\frac{(\hat {\textbf x}_1^T\textbf l_2)(\hat {\textbf x}_3^T\textbf l_2)}{(\hat {\textbf x}_1^T\textbf l_1)(\hat {\textbf x}_3^T\textbf l_1)}, \quad \quad \lambda_3=\frac{(\hat {\textbf x}_2^T\textbf l_2)(\hat {\textbf x}_3^T\textbf l_2)}{(\hat {\textbf x}_2^T\textbf l_1)(\hat {\textbf x}_3^T\textbf l_1)},
\end{aligned}
\end{equation}
are three distinct real eigenvalues of ${\textbf C_{s,t}^{(1)}}^{-1}\textbf C_{s,t}^{(2)}$ (See Appendix A). Further, we may assert that $\textbf C_{s,t}^{(1)}$ and $\textbf C_{s,t}^{(2)}$ must have four distinct intersection points. That means that the problem of determining a non-degenerate conic must have four distinct solutions. This completes the proof of the theorem. \hfill {$\Box$}\\

\noindent\textbf {Proposition 3.} \emph{Given three points} $\hat {\textbf x}_{1}$, $\hat {\textbf x}_{2}$, $\hat {\textbf x}_{3}\in\mathrm P\mathbb{R}^2$ \emph{and two lines} $\textbf l_1$, $\textbf l_2\in\mathrm P\mathbb{R}^2$ \emph{in general position such that} $\hat {\textbf x}_i^T\textbf l_1\neq0$, $\hat {\textbf x}_i^T\textbf l_2\neq0$, $i=1,2,3$ \emph{and the point} $\textbf p=\textbf l_1\times\textbf l_2$ \emph{does not lie on any side of the triangle} $\hat {\textbf x}_1\hat {\textbf x}_2\hat {\textbf x}_3$. \emph{Then the problem of determining a non-degenerate conic must have four distinct real solutions if and only if the signs of} $(\hat {\textbf x}_1^T\textbf l_1)(\hat {\textbf x}_1^T\textbf l_2)$, $(\hat {\textbf x}_2^T\textbf l_1)(\hat {\textbf x}_2^T\textbf l_2)$, $(\hat {\textbf x}_3^T\textbf l_1)(\hat {\textbf x}_3^T\textbf l_2)$ \emph{are same}.\\

\noindent\textbf {Proof.} On one hand, from (15) in Theorem 4, we know that $\lambda_i\textbf C_{s,t}^{(1)}-\textbf C_{s,t}^{(2)}$, $i=1,2,3$ are three line-pairs such that the intersection of any two of them gives the intersections of $\textbf C_{s,t}^{(1)}$, $\textbf C_{s,t}^{(2)}$.

On the other hand, we may check
\begin{equation*}
\begin{aligned}
\lambda_i\textbf C_{s,t}^{(1)}-\textbf C_{s,t}^{(2)}=\textbf M^T(\lambda_i{\textbf C'}_{s,t}^{(1)}-{\textbf C'}_{s,t}^{(2)})\textbf M,
\end{aligned}
\end{equation*}
where
\begin{equation*}
\begin{aligned}
&\textbf M=\left[\begin{array}{ccc}
    det([\hat {\textbf x}_1, \textbf p, \hat {\textbf x}_3]) & \frac{det([\hat {\textbf x}_1, \textbf p, \hat {\textbf x}_3][\hat {\textbf x}_{123}])}{det([\textbf p, \hat {\textbf x}_2, \hat {\textbf x}_3])} & 0\\  
    -det([\textbf p, \hat {\textbf x}_2, \hat {\textbf x}_3]) & det([\hat {\textbf x}_{123}]) & 0 \\  
    0 & 0 & det([\hat {\textbf x}_1, \hat {\textbf x}_2, \textbf p]) \\  
  \end{array}
\right],
\end{aligned}
\end{equation*}
and
\begin{equation*}
\begin{aligned}
&{\textbf C'}_{s,t}^{(i)}=(\hat {\textbf x}_1^T\textbf l_i)\left[\begin{array}{ccc}
    -k(\hat {\textbf x}_2^T\textbf l_i) & 0 & 0\\  
    0 & \hat {\textbf x}_1^T\textbf l_i+\frac{\hat {\textbf x}_2^T\textbf l_i}{k} & -\hat {\textbf x}_3^T\textbf l_i \\  
    0 & -\hat {\textbf x}_3^T\textbf l_i & \frac{(\hat {\textbf x}_3^T\textbf l_i)^2}{\hat {\textbf x}_1^T\textbf l_i} \\  
  \end{array}
\right], i=1,2,\\
&k=\frac{det([\textbf p, \hat {\textbf x}_2, \hat {\textbf x}_3])}{det([\hat {\textbf x}_1, \textbf p, \hat {\textbf x}_3])}.
\end{aligned}
\end{equation*}
This tells us that for $i=1,2,3$ the matrix $\lambda_i\textbf C_{s,t}^{(1)}-\textbf C_{s,t}^{(2)}$ is congruent to the matrix $\lambda_i{\textbf C'}_{s,t}^{(1)}-{\textbf C'}_{s,t}^{(2)}$. Consequently, they have the same positive and negative inertial indexes according to inertial theorem [12]. By careful computation, we may know that\\

\par\setlength\parindent{0.5em}(i) the product of two non-zero eigenvalues of the degenerate conic $\lambda_1{\textbf C'}_{s,t}^{(1)}-{\textbf C'}_{s,t}^{(2)}$ is as follows:\\
\begin{equation*}
\begin{aligned}
-\lambda_1((\hat {\textbf x}_1^T\textbf l_1)(\hat {\textbf x}_3^T\textbf l_2)-(\hat {\textbf x}_1^T\textbf l_2)(\hat {\textbf x}_3^T\textbf l_1))^2,\\
\end{aligned}
\end{equation*}

\par\setlength\parindent{0.2em}(ii) the product of two non-zero eigenvalues of the degenerate conic $\lambda_2{\textbf C'}_{s,t}^{(1)}-{\textbf C'}_{s,t}^{(2)}$ is as follows:\\
\begin{equation*}
\begin{aligned}
-\lambda_2k^2((\hat {\textbf x}_1^T\textbf l_1)(\hat {\textbf x}_3^T\textbf l_2)-(\hat {\textbf x}_1^T\textbf l_2)(\hat {\textbf x}_3^T\textbf l_1))^2,\\
\end{aligned}
\end{equation*}

(iii) the product of two non-zero eigenvalues of the degenerate conic $\lambda_3{\textbf C'}_{s,t}^{(1)}-{\textbf C'}_{s,t}^{(2)}$ is as follows:\\
\begin{equation*}
\begin{aligned}
-\lambda_3k^2&((\hat {\textbf x}_1^T\textbf l_1)(\hat {\textbf x}_3^T\textbf l_2)-(\hat {\textbf x}_1^T\textbf l_2)(\hat {\textbf x}_3^T\textbf l_1))^2(\frac{((\hat {\textbf x}_1^T\textbf l_1)(\hat {\textbf x}_2^T\textbf l_2)-(\hat {\textbf x}_1^T\textbf l_2)(\hat {\textbf x}_2^T\textbf l_1))^2}{((\hat {\textbf x}_3^T\textbf l_1)(\hat {\textbf x}_2^T\textbf l_2)-(\hat {\textbf x}_3^T\textbf l_2)(\hat {\textbf x}_2^T\textbf l_1))^2}+1).
\end{aligned}
\end{equation*}

According to (15), we may know that the signatures of $\lambda_i{\textbf C'}_{s,t}^{(1)}-{\textbf C'}_{s,t}^{(2)}$, $i=1,2,3$ with the rank 2 are zero if and only if the signs of $(\hat {\textbf x}_1^T\textbf l_1)(\hat {\textbf x}_1^T\textbf l_2)$, $(\hat {\textbf x}_2^T\textbf l_1)(\hat {\textbf x}_2^T\textbf l_2)$, $(\hat {\textbf x}_3^T\textbf l_1)(\hat {\textbf x}_3^T\textbf l_2)$ are same. Further, according to the inertial theorem, we may assert that the signatures of $\lambda_i{\textbf C}_{s,t}^{(1)}-{\textbf C}_{s,t}^{(2)}$, $i=1,2,3$ with the rank 2 are also zero, that means that $\lambda_i{\textbf C}_{s,t}^{(1)}-{\textbf C}_{s,t}^{(2)}$, $i=1,2,3$ are three real line-pairs, and further the intersection of any two of them will give four distinct real intersections of ${\textbf C}_{s,t}^{(1)}$, ${\textbf C}_{s,t}^{(2)}$. This completes the proof of the proposition. \hfill {$\Box$}\\

In the following, we will study all sorts of special configurations of three points and two lines in general position.\\

\noindent\textbf {Theorem 5.} \emph{Given three points and two lines in general position in} $\mathrm P\mathbb{R}^2$. \emph{If there exist some two points lying on the two lines respectively, then the problem of determining a non-degenerate conic has a unique real solution}.\\

\noindent\textbf {Proof.} In this case, we allocate $\hat {\textbf x}_2^T\textbf l_1=0$, $\hat {\textbf x}_3^T\textbf l_2=0$, $\hat {\textbf x}_1^T\textbf l_i\neq0$, $i=1,2$ and let $\textbf x_4=t\hat {\textbf x}_1+\textbf p$, where $\textbf p=\textbf l_1\times\textbf l_2$, then we have
\begin{equation*}
\begin{aligned}
&det([\textbf p, \hat {\textbf x}_2, \hat {\textbf x}_3]\hat {\textbf x}^T_1\textbf l_1+det([\hat {\textbf x}_1, \hat {\textbf x}_2, \textbf p]\hat {\textbf x}^T_3\textbf l_1=0,\\
&det([\textbf p, \hat {\textbf x}_2, \hat {\textbf x}_3]\hat {\textbf x}^T_1\textbf l_2+det([\hat {\textbf x}_1, \textbf p, \hat {\textbf x}_2]\hat {\textbf x}^T_2\textbf l_2=0.\\
\end{aligned}
\end{equation*}
Then according to (12) and (13), we may get a unique intersection point of the conics $\textbf C^{(1)}_{s,t}$ and $\textbf C^{(2)}_{s,t}$ as follows:
\begin{equation}
\begin{aligned}
s=\frac{1}{2}, \quad \quad t=-\frac{det([\textbf p, \hat {\textbf x}_2, \hat {\textbf x}_3])}{2det([\hat {\textbf x}_{123}])}.
\end{aligned}
\end{equation}
So the problem of determining a non-degenerate conic has a unique real solution. This completes the proof of the theorem. \hfill {$\Box$}\\

\noindent\textbf {Theorem 6.} \emph{Given three points and two lines in general position in} $\mathrm P\mathbb{R}^2$. \emph{If there exist some two points and the intersection point of two lines being collinear, then}\\

\par\setlength\parindent{0.5em}(i) \emph{when the remaining third point does not lie on any of the two lines, the problem of determining a non-degenerate conic must have two distinct solutions}.\\

\par\setlength\parindent{0.2em}(ii) \emph{when the remaining third point lies on some one of the two lines, the problem of determining a non-degenerate conic has a unique real solution}.\\

\noindent\textbf {Proof.} For (i), in this case, we allocate $\hat {\textbf x}^T_1\textbf l_1\neq0$,$\hat {\textbf x}^T_1\textbf l_2\neq0$ and $\hat {\textbf x}_2$, $\hat {\textbf x}_3$, $\textbf p$ be collinear, where $\textbf p=\textbf l_1\times\textbf l_2$.

Let $\textbf x_4=t\hat {\textbf x}_1+\textbf p$. Then we may deduce that $\boldsymbol {\xi}_{3}^T\textbf l_1=\boldsymbol {\xi}_{3}^T\textbf l_2=0$. According to (11), we have
\begin{equation}
\begin{aligned}
s=\frac{(\boldsymbol {\xi}_{1}^T\textbf l_1)^2}{(\boldsymbol {\xi}_{1}^T\textbf l_1)^2-(\boldsymbol {\xi}_{2}^T\textbf l_1)^2}=\frac{(\boldsymbol {\xi}_{1}^T\textbf l_2)^2}{(\boldsymbol {\xi}_{1}^T\textbf l_2)^2-(\boldsymbol {\xi}_{2}^T\textbf l_2)^2}.
\end{aligned}
\end{equation}
Further, from (1) we may get
\begin{equation}
\begin{aligned}
t^2=\frac{det([\hat {\textbf x}_1, \hat {\textbf x}_2, \textbf p])^2(\hat {\textbf x}_3^T\textbf l_1)(\hat {\textbf x}_3^T\textbf l_2)}{det([\hat {\textbf x}_{123}])^2(\hat {\textbf x}_1^T\textbf l_1)(\hat {\textbf x}_1^T\textbf l_2)}.
\end{aligned}
\end{equation}
Since the right of (18) must be not zero, $t$ has two distinct roots and further we may get two corresponding $s$ from (1) and (17). So the problem of determining a non-degenerate conic has two distinct solutions.\\

For (ii), in this case, we allocate $\hat {\textbf x}_1$, $\hat {\textbf x}_2$, $\textbf p$ be collinear and $\hat {\textbf x}^T_3\textbf l_1=0$, where $\textbf p=\textbf l_1\times\textbf l_2$. It is worth noting that we can not let $\textbf x_4=t\hat {\textbf x}_1+\textbf p$, but let $\textbf x_4=t\textbf p+\textbf q$, where $\textbf q=\textbf p\times\textbf l_2$. Now since $\hat {\textbf x}_3^T\textbf l_1=0$ and $\textbf x_4^T\textbf l_2=0$, according to Proposition 1, from (10) and (1) we may solve
\begin{equation}
\begin{aligned}
&s=2,\quad \quad t=\frac{det([\textbf q, \hat {\textbf x}_2, \hat {\textbf x}_3])(\hat {\textbf x}_1^T\textbf l_2)-det([\hat {\textbf x}_1, \textbf q, \hat {\textbf x}_3])(\hat {\textbf x}_2^T\textbf l_2)}{2det([\hat {\textbf x}_1, \textbf p, \hat {\textbf x}_3])(\hat {\textbf x}_2^T\textbf l_2)}.
\end{aligned}
\end{equation}
So the problem of determining a non-degenerate conic has a unique real solution. This completes the proof of the theorem. \hfill {$\Box$}\\

According to the proof of (i) in Theorem 6, we may deduce immediately that under what configuration the problem of determining a non-degenerate conic has two distinct real solutions. More formally, we have the following proposition.\\

\noindent\textbf {Proposition 4.} \emph{Given three points} $\hat {\textbf x}_{1}$, $\hat {\textbf x}_{2}$, $\hat {\textbf x}_{3}\in\mathrm P\mathbb{R}^2$ \emph{and two lines} $\textbf l_1$, $\textbf l_2\in\mathrm P\mathbb{R}^2$ \emph{in general position such that} $\hat {\textbf x}_1^T\textbf l_i\neq0$, $i=1,2$ \emph{and} $\hat {\textbf x}_2$, $\hat {\textbf x}_3$, $\textbf l_1\times\textbf l_2$ \emph{are collinear. Then the problem of determining a non-degenerate conic has two distinct real solutions if and only if} $(\hat {\textbf x}_1^T\textbf l_1)(\hat {\textbf x}_3^T\textbf l_1)$ \emph{and} $(\hat {\textbf x}_1^T\textbf l_2)(\hat {\textbf x}_3^T\textbf l_2)$ \emph{have same signs}.\\

\noindent\textbf {Theorem 7.} \emph{Given three points and two lines in general position in} $\mathrm P\mathbb{R}^2$. \emph{If there exists only one point lying on some line and the remaining two points and the intersection point of two lines are not collinear, then the problem of determining a non-degenerate conic has two distinct solutions}.\\

\noindent\textbf {Proof.} In this case, we allocate $\hat {\textbf x}_3^T\textbf l_1=0$ and let $\textbf x_4=t\hat {\textbf x}_1+\textbf p$, where $\textbf p=\textbf l_1\times\textbf l_2$, then according to Proposition 1, from (10) and (1), we may have
\begin{equation}
\begin{aligned}
s=\frac{(\boldsymbol {\xi}_{1}^T\textbf l_1)^2-(\boldsymbol {\xi}_{2}^T\textbf l_1)^2}{2(\boldsymbol {\xi}_{3}^T\textbf l_1)^2}+\frac{1}{2}=-\frac{det([\hat {\textbf x}_{123}])}{det([\textbf p, \hat {\textbf x}_{2}, \hat {\textbf x}_{3}])}t.
\end{aligned}
\end{equation}
After we substitute it into the second equation in (12), we may derive a quadratic equation about $t$ as follows:
\begin{equation}
\begin{aligned}
&4det([\hat {\textbf x}_{123}])^2(\hat {\textbf x}_{1}^T\textbf l_2)t^2+4det([\hat {\textbf x}_{123}])det([\textbf p, \hat {\textbf x}_{2}, \hat {\textbf x}_{3}])(\hat {\textbf x}_{1}^T\textbf l_2)t\\
&\quad \quad \quad \quad \quad \quad \quad \quad \quad \quad \quad -det([\textbf p, \hat {\textbf x}_{2}, \hat {\textbf x}_{3}])det([\hat {\textbf x}_{1}, \hat {\textbf x}_{2}, \textbf p])(\hat {\textbf x}_{3}^T\textbf l_2)=0.
\end{aligned}
\end{equation}
Further, by the detailed derivation, the discriminant of (21) may be represented as follows:
\begin{equation}
\begin{aligned}
&\Delta=-16det([\hat {\textbf x}_{123}])^2det([\textbf p, \hat {\textbf x}_{2}, \hat {\textbf x}_{3}])det([\hat {\textbf x}_{1}, \textbf p, \hat {\textbf x}_{3}])(\hat {\textbf x}_{1}^T\textbf l_2)(\hat {\textbf x}_{2}^T\textbf l_2).
\end{aligned}
\end{equation}
Since the discriminant (22) must be nonzero, we may assert that the problem of determining a non-degenerate conic must have two distinct solutions. This completes the proof of the theorem. \hfill {$\Box$}\\

An immediate consequence of the Theorem 7 is that the condition of two distinct real solutions can be determined. More formally, we have the following proposition.\\

\noindent\textbf {Proposition 5.} \emph{Given three points} $\hat {\textbf x}_{1}$, $\hat {\textbf x}_{2}$, $\hat {\textbf x}_{3}\in\mathrm P\mathbb{R}^2$ \emph{and two lines} $\textbf l_1$, $\textbf l_2\in\mathrm P\mathbb{R}^2$ \emph{in general position such that} $\hat {\textbf x}_3^T\textbf l_1=0$ \emph{and} $\hat {\textbf x}_1$, $\hat {\textbf x}_2$, $\textbf p=\textbf l_1\times\textbf l_2$ \emph{are not collinear. Then the problem of determining a non-degenerate conic has two distinct real solutions if and only if}
\begin{equation*}
\begin{aligned}
det([\textbf p, \hat {\textbf x}_{2}, \hat {\textbf x}_{3}])det([\hat {\textbf x}_{1}, \textbf p, \hat {\textbf x}_{3}])(\hat {\textbf x}_{1}^T\textbf l_2)(\hat {\textbf x}_{2}^T\textbf l_2)<0.\\
\end{aligned}
\end{equation*}

\subsection{The Rest of All Minimal Configurations}

Since the duality principle of the projective plane, we may obtain easily the closed-form solutions to the non-degenerate dual conic from the cases of two-point + three-line configuration, one-point + four-line configuration and five-line configuration.

\section{Implementation}

In this section, we summarize the proposed algorithms in the following.\\

\noindent\textbf {Algorithm 1: conic from five points}\\

\par\setlength\parindent{0.5em}(i) Given five points in general position. Denote them arbitrarily by $\hat {\textbf x}_i$, $i=1,...,5$. Then, compute $\boldsymbol {\xi}_{1}$, $\boldsymbol {\xi}_{2}$, $\boldsymbol {\xi}_{3}$ according to (1).

\par\setlength\parindent{0.2em}(ii) Compute the conic $\textbf C$ according to (3).\\

\noindent\textbf {Algorithm 2: conics from four points \!+\! one line}\\

\noindent\textbf {Case 1:} Given four points and one line in general position such that there exists one point lying on the line.\\

\par\setlength\parindent{0.5em}(i) Denote four points arbitrarily by $\hat {\textbf x}_i$, $i=1,...,4$ and denote the line by $\textbf l$. Then, compute $\boldsymbol {\xi}_{1}$, $\boldsymbol {\xi}_{2}$, $\boldsymbol {\xi}_{3}$ according to (1).

\par\setlength\parindent{0.2em}(ii) Compute $s$ according to (10) and further compute the conic $\textbf C$ according to (2).\\

\noindent\textbf {Case 2:} Given four points and one line in general position such that the line passes through some one of three vertices of the diagonal triangle of the quadrangle consisting of given four points.\\

\par\setlength\parindent{0.5em}(i) Denote the line by $\textbf l$ and denote four points by $\hat {\textbf x}_i$, $i=1,...,4$ such that $(\hat {\textbf x}_1\times\hat {\textbf x}_4)\times(\hat {\textbf x}_2\times\hat {\textbf x}_3)$ lies on $\textbf l$. Then, compute $\boldsymbol {\xi}_{1}$, $\boldsymbol {\xi}_{2}$, $\boldsymbol {\xi}_{3}$ according to (1).

\par\setlength\parindent{0.2em}(ii) Compute $s$ according to (11) and further compute the conic $\textbf C$ according to (2).\\

\noindent\textbf {Case 3:}  Given four points and one line in general position such that the line neither passes through any one of four points and nor passes through any one of three vertices of the diagonal triangle of the quadrangle consisting of given four points.\\

\par\setlength\parindent{0.5em}(i) Denote the line $\textbf l$ and denote four points arbitrarily by $\hat {\textbf x}_i$, $i=1,...,4$. Then, compute $\boldsymbol {\xi}_{1}$, $\boldsymbol {\xi}_{2}$, $\boldsymbol {\xi}_{3}$ according to (1).

\par\setlength\parindent{0.2em}(ii) Compute $s_1$ and $s_2$ according to (8).

(iii) Compute the conics $\textbf C_1$ and $\textbf C_2$ according to (7).\\

\noindent\textbf {Algorithm 3: conics from three points\! +\! two lines}\\

\noindent\textbf {Case 1:} Given three points and two lines in general position such that there exist two points lying on two lines respectively.\\

\par\setlength\parindent{0.5em}(i) Denote the point not lying on $\textbf l_1$ or $\textbf l_2$ by $\hat {\textbf x}_1$ and denote two lines arbitrarily by $\textbf l_1$ and $\textbf l_2$. Further, denote the point lying on $\textbf l_1$ by $\hat {\textbf x}_2$ and denote the point lying on $\textbf l_2$ by $\hat {\textbf x}_3$.

\par\setlength\parindent{0.2em}(ii) Compute $\textbf p=\textbf l_1\times\textbf l_2$ and further compute $t$ according to (16). Then compute $\boldsymbol {\xi}_1$, $\boldsymbol {\xi}_2$ and $\boldsymbol {\xi}_3$ according to $\textbf x_4=t\hat {\textbf x}_1+\textbf p$ and (1).

(iii) Let $s=\frac{1}{2}$ and further compute the conic $\textbf C$ according to (2).\\

\noindent\textbf {Case 2:} Given three points and two lines in general position such that some two points and the intersection points of two lines are collinear and the remaining third point lies on some one of two lines.\\

\par\setlength\parindent{0.5em}(i) Denote two lines and one point by $\textbf l_1$, $\textbf l_2$ and $\hat {\textbf x}_3$ such that $\hat {\textbf x}_3^T\textbf l_1=0$. Then denote the other two points arbitrarily by $\hat {\textbf x}_1$ and $\hat {\textbf x}_2$.

\par\setlength\parindent{0.2em}(ii) Compute $\textbf p=\textbf l_1\times\textbf l_2$, $\textbf q=\textbf p\times\textbf l_2$ and further compute $t$ according to (19). Then compute $\boldsymbol {\xi}_1$, $\boldsymbol {\xi}_2$ and $\boldsymbol {\xi}_3$ according to $\textbf x_4=t\textbf p+\textbf q$ and (1).

(iii) Let $s=2$ and further compute the conic $\textbf C$ according to (2).\\

\noindent\textbf {Case 3:} Given three points and two lines in general position such that some two points and the intersection points of two lines are collinear and the remaining third point does not lie on any one of two lines.\\

\par\setlength\parindent{0.5em}(i) Denote two lines arbitrarily by $\textbf l_1$, $\textbf l_2$ and denote two points that are collinear with the intersection point $\textbf l_1\times\textbf l_2$ arbitrarily by $\hat {\textbf x}_2$ and $\hat {\textbf x}_3$. Then denote the remaining third point by $\hat {\textbf x}_1$.

\par\setlength\parindent{0.2em}(ii) Compute $\textbf p=\textbf l_1\times\textbf l_2$ and further compute $t_1$ and $t_2$ according to (18). Then compute $\boldsymbol {\xi}^{(i)}_1$, $\boldsymbol {\xi}^{(i)}_2$, $\boldsymbol {\xi}^{(i)}_3$, $i=1,2$ according to $\textbf x^{(i)}_4=t_i\hat {\textbf x}_1+\textbf p$, $i=1,2$ and (1).

(iii) Compute $s_1$ and $s_2$ according to (17) and further compute the conics $\textbf C_1$ and $\textbf C_2$ according to (2).\\

\noindent\textbf {Case 4:} Given three points and two lines in general position such that some one point lies on some one line and the remaining two points and the intersection point of two lines are not collinear.\\

\par\setlength\parindent{0.5em}(i) Denote two points not lying on the lines arbitrarily by $\hat {\textbf x}_1$ and $\hat {\textbf x}_2$ and denote the point lying on some one line by $\hat {\textbf x}_3$. Denote the line through $\hat {\textbf x}_3$ by $\textbf l_1$ and another line by $\textbf l_2$.

\par\setlength\parindent{0.2em}(ii) Compute $\textbf p=\textbf l_1\times\textbf l_2$ and further compute $t_1$ and $t_2$ according to (21). Then compute $\boldsymbol {\xi}^{(i)}_1$, $\boldsymbol {\xi}^{(i)}_2$, $\boldsymbol {\xi}^{(i)}_3$, $i=1,2$ according to $\textbf x^{(i)}_4=t_i\hat {\textbf x}_1+\textbf p$, $i=1,2$ and (1).

(iii) Compute $s_1$ and $s_2$ according to (20) and further compute the conics $\textbf C_1$ and $\textbf C_2$ according to (2).\\

\noindent\textbf {Case 5:} Given three points and two lines in general position such that all the points do not lie on the lines and furthermore the line through any two points does not go through the intersection point of two lines.\\

\par\setlength\parindent{0.5em}(i) Denote two lines arbitrarily by $\textbf l_1$, $\textbf l_2$ and denote three points arbitrarily by $\hat {\textbf x}_1$, $\hat {\textbf x}_2$, $\hat {\textbf x}_3$.

\par\setlength\parindent{0.2em}(ii) Compute $\textbf p=\textbf l_1\times\textbf l_2$ and further compute the intersection points $(s_i, t_i,1)^T$, $i=1,2,3,4$ of two non-degenerate conics in (12). Then compute $\boldsymbol {\xi}^{(i)}_1$, $\boldsymbol {\xi}^{(i)}_2$, $\boldsymbol {\xi}^{(i)}_3$, $i=1,2,3,4$ according to $\textbf x^{(i)}_4=t_i\hat {\textbf x}_1+\textbf p$, $i=1,2,3,4$ and (1).

(iii) Compute the conics $\textbf C_i$, $i=1,2,3,4$ according to (2).\\

\section{Experiments}

In this section, we show some examples of determining the conics from the point-line minimal configurations to demonstrate the correctness and advantage of our results.

Fist, given five points $\hat {\textbf x}_i$, $i=1,2,3,4,5$ in general position, we may obtain directly a real non-degenerate conic conic by Theorem 2 and Algorithm 1, as shown in Fig. 2, where Fig. 2(a) shows an ellipse, Fig. 2(b) shows a hyperbola and Fig. 2(c) shows a parabola.
\begin{figure}
\begin{center}
\subfigure[]{
\label{fig:subfig:a}  
\includegraphics[width=0.3\textwidth]{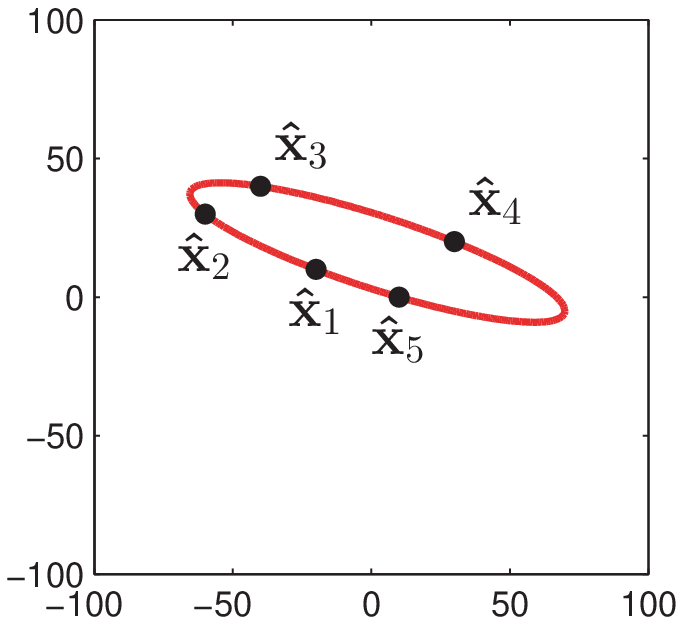}}
\subfigure[]{
\label{fig:subfig:b}  
\includegraphics[width=0.3\textwidth]{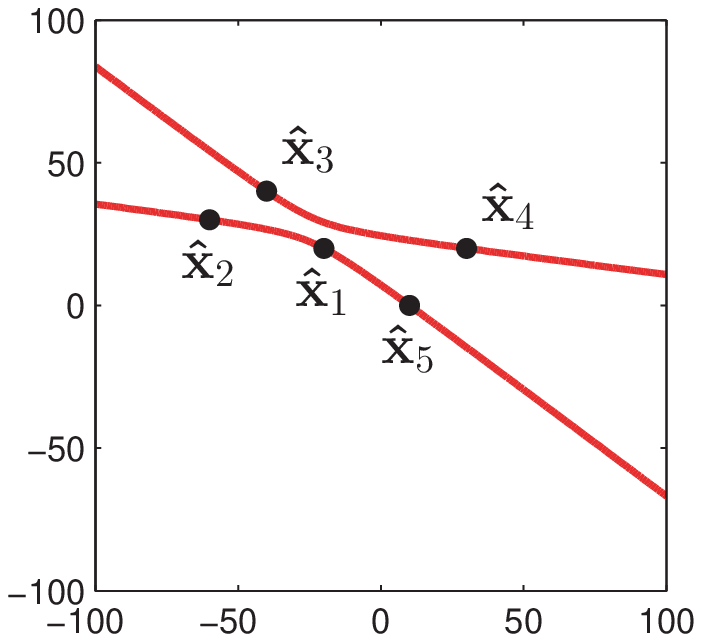}}
\subfigure[]{
\label{fig:subfig:c}  
\includegraphics[width=0.3\textwidth]{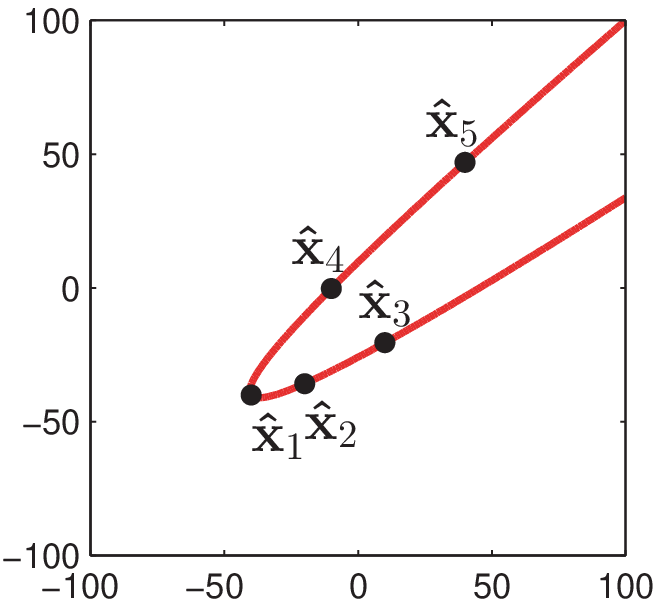}}
\caption{Conics from five points $\hat {\textbf x}_i$, $i=1,2,3,4,5$ in general position. (a), (b) and (c) show a unique real non-degenerate conic respectively according to Theorem 2 and Algorithm 1.}
\label{fig:subfig}  
\end{center}
\end{figure}

Second, given four points $\hat {\textbf x}_i$, $i=1,2,3,4$ and one line $\textbf l$ in general position, according to Proposition 1 and Proposition 2, we may determine easily whether there exists at least a real non-degenerate conic passing through $\hat {\textbf x}_i$, $i=1,2,3,4$ and tangent to $\textbf l$. If yes, then we may obtain one or two real non-degenerate conics by Theorem 3 and Algorithm 2, as shown in Fig. 3, where Fig. 3(a) shows a unique real non-degenerate conic by Proposition 1 (Note that the line $\textbf l$ passes through the point $\hat {\textbf x}_1$.) and Case 1 of Algorithm 2. Fig. 3(b) shows a unique real non-degenerate conic by Proposition 1 (Note that the line $\textbf l$ passes through the vertex $\boldsymbol {\xi}_{1}$ of the diagonal triangle of the quadrangle $\hat {\textbf x}_1\hat {\textbf x}_2\hat {\textbf x}_3\hat {\textbf x}_4$.) and Case 2 of Algorithm 2. In Fig. 3(c), all four points lie on the same side of the line, but $det([\hat {\textbf x}_{123}])$, $det([\hat {\textbf x}_{124}])$, $det([\hat {\textbf x}_{134}])$ have same directions and $det([\hat {\textbf x}_{234}])$ is anti-direction, so there exist no real solutions by Proposition 2. Fig. 3(d), Fig. 3(e) and Fig. 3(f) show two distinct real non-degenerate conics by Proposition 2 and Case 3 of Algorithm 2.
\begin{figure}
\begin{center}
\subfigure[]{
\label{fig:subfig:a}  
\includegraphics[width=0.3\textwidth]{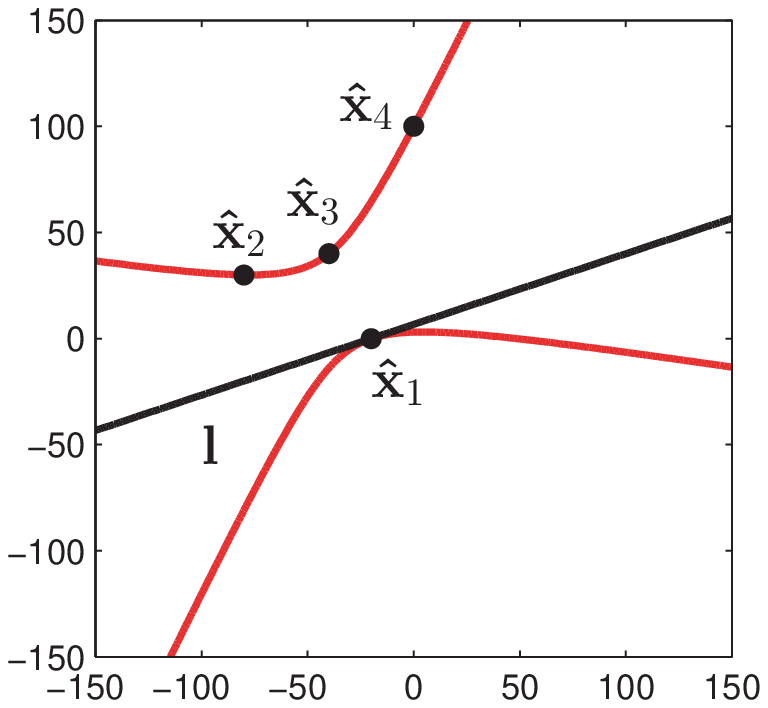}}
\subfigure[]{
\label{fig:subfig:b}  
\includegraphics[width=0.3\textwidth]{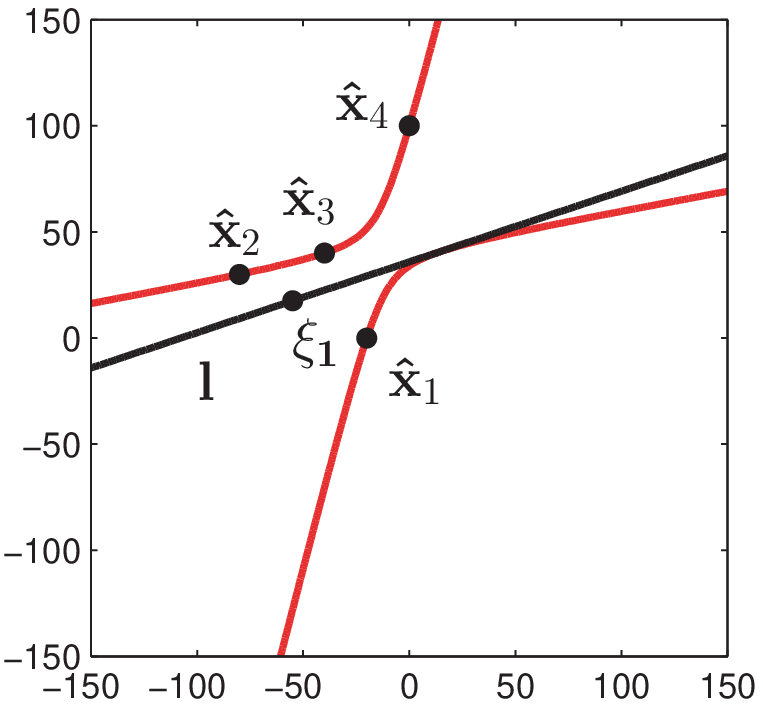}}
\subfigure[]{
\label{fig:subfig:c}  
\includegraphics[width=0.3\textwidth]{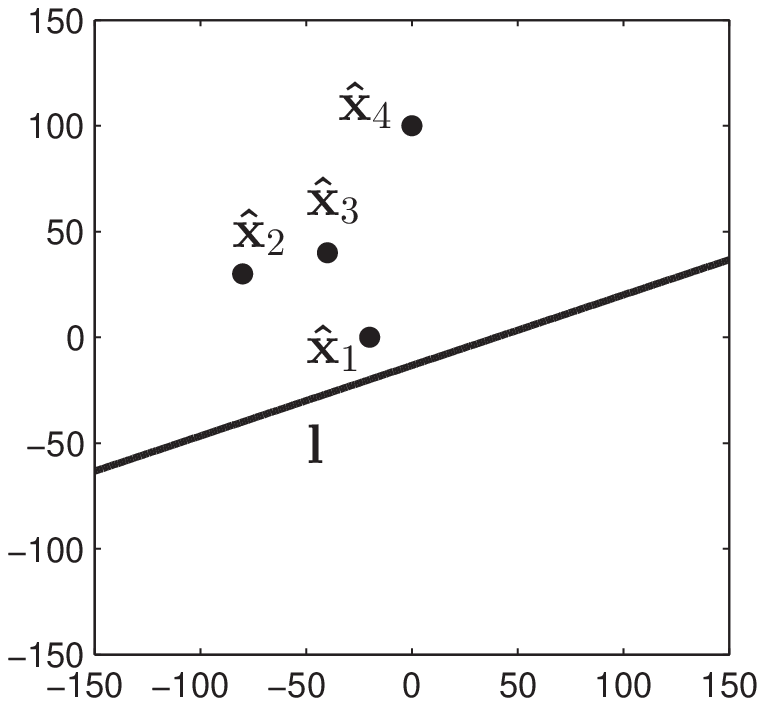}}
\vspace{1mm}
\subfigure[]{
\label{fig:subfig:a}  
\includegraphics[width=0.3\textwidth]{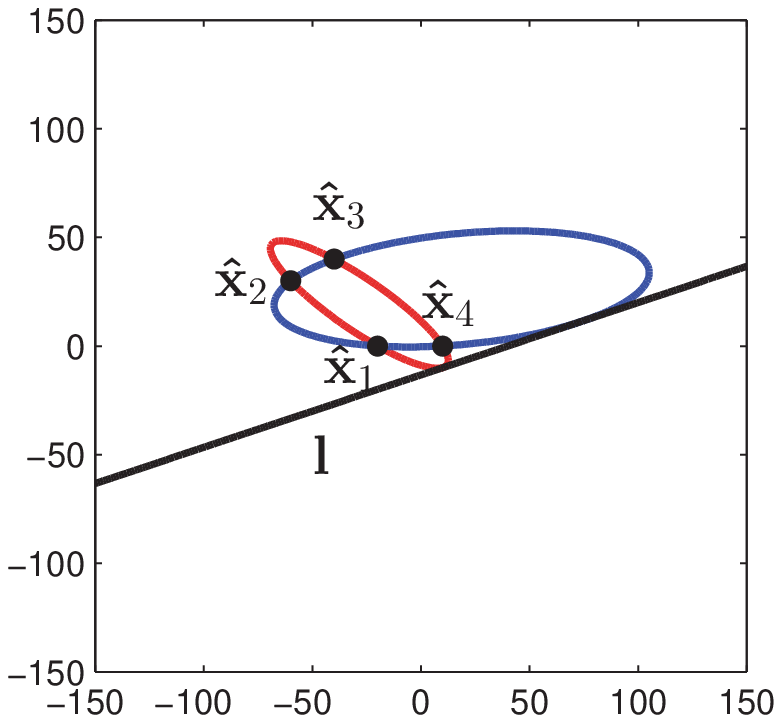}}
\subfigure[]{
\label{fig:subfig:b}  
\includegraphics[width=0.3\textwidth]{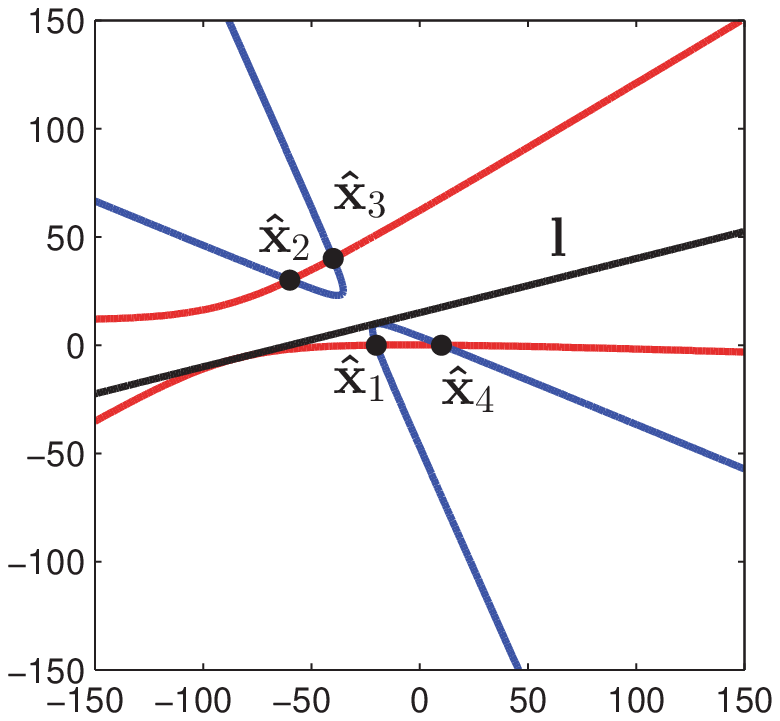}}
\subfigure[]{
\label{fig:subfig:c}  
\includegraphics[width=0.3\textwidth]{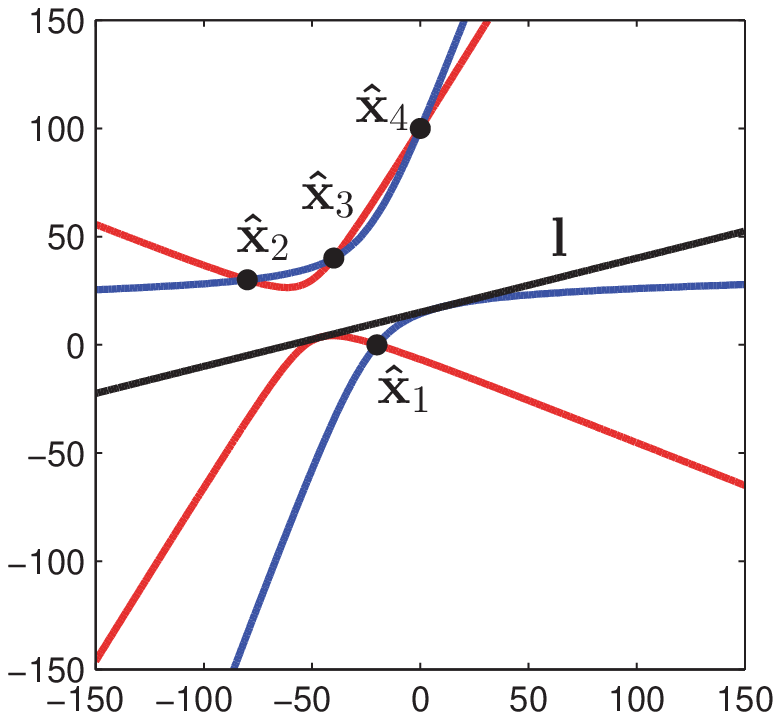}}
\caption{Conics from four points $\hat {\textbf x}_1$, $\hat {\textbf x}_2$, $\hat {\textbf x}_3$, $\hat {\textbf x}_4$ and one line $\textbf l$ in general position. (a) shows a unique real non-degenerate conic by Proposition 1 and Case 1 of Algorithm 2, where $\textbf l$ passes through $\hat {\textbf x}_1$. (b) shows a unique real non-degenerate conic by Proposition 1 and Case 2 of Algorithm 2, where $\textbf l$ passes through the vertex $\boldsymbol {\xi}_{1}$ of the diagonal triangle of the quadrangle $\hat {\textbf x}_1\hat {\textbf x}_2\hat {\textbf x}_3\hat {\textbf x}_4$. (c) shows no real solutions by Proposition 2. (d), (e) and (f) show two distinct real non-degenerate conics by Proposition 2 and Case 3 of Algorithm 2.}
\label{fig:subfig}  
\end{center}
\end{figure}

Last, given three points $\hat {\textbf x}_i$, $i=1,2,3$ and two lines $\textbf l_1$, $\textbf l_2$ in general position, by Proposition 3, Theorem 5, Theorem 6, Proposition 4, Proposition 5, we may determine easily whether there exists at least a real non-degenerate conic passing through $\hat {\textbf x}_i$, $i=1,2,3$ and tangent to $\textbf l_1$, $\textbf l_2$. If yes, then we may obtain one, two or four real non-degenerate conics by Algorithm 3, as shown in Fig. 4, where Fig. 4(a) shows a unique solution by Theorem 5 and Case 1 of Algorithm 3. Fig. 4(b) shows a unique solution by (ii) of Theorem 6 and Case 2 of Algorithm 3. Fig. 4(c) and Fig. 4(d) show two distinct real solutions by (i) of Theorem 6 and Case 3 of Algorithm 3. In Fig. 4(e), because the points $\hat {\textbf x}_1$, $\hat {\textbf x}_3$ lie on the same side of the line $\textbf l_1$, but they lie on the opposite sides of the line $\textbf l_2$, that means that the signs of $(\hat {\textbf x}^T_1\textbf l_1)(\hat {\textbf x}^T_3\textbf l_1)$ and $(\hat {\textbf x}^T_1\textbf l_2)(\hat {\textbf x}^T_3\textbf l_2)$ are different, so there exist no real solutions by Proposition 4. Fig. 4(f) and Fig. 4(g) show two distinct real solutions by Proposition 5 and Case 4 of Algorithm 3. In Fig. 4(h), on one hand, because the points $\hat {\textbf x}_1$, $\hat {\textbf x}_2$ lie on the opposite sides of the line $\textbf l_2$, that means that $(\hat {\textbf x}^T_1\textbf l_2)(\hat {\textbf x}^T_2\textbf l_2)<0$. On the other hand, the directions of $det([\textbf l_1\times\textbf l_2, \hat {\textbf x}_{2}, \hat {\textbf x}_{3}])$ and $det([\hat {\textbf x}_{1}, \textbf l_1\times\textbf l_2, \hat {\textbf x}_{3}])$ are different, that means that
\begin{equation*}
\begin{aligned}
det([\textbf l_1\times\textbf l_2, \hat {\textbf x}_{2}, \hat {\textbf x}_{3}])det([\hat {\textbf x}_{1}, \textbf l_1\times\textbf l_2, \hat {\textbf x}_{3}])<0,
\end{aligned}
\end{equation*}
so there exist no real solutions by Proposition 5. Fig. 4(i) and Fig. 4(j) show four distinct real solutions by Proposition 3 and Case 5 of Algorithm 3. In Fig. 4(k) and Fig. 4(l), both the point $\hat {\textbf x}_2$ and the point $\hat {\textbf x}_3$ lie on the same side of the lines $\textbf l_1$, $\textbf l_2$, but the point $\hat {\textbf x}_1$ does not. That means that the signs of $(\hat {\textbf x}^T_1\textbf l_1)(\hat {\textbf x}^T_1\textbf l_2)$, $(\hat {\textbf x}^T_2\textbf l_1)(\hat {\textbf x}^T_2\textbf l_2)$, $(\hat {\textbf x}^T_3\textbf l_1)(\hat {\textbf x}^T_3\textbf l_2)$ are different, so there exist no real solutions by Proposition 3.

\begin{figure}
\begin{center}
\subfigure[]{
\label{fig:subfig:a}  
\includegraphics[width=0.26\textwidth]{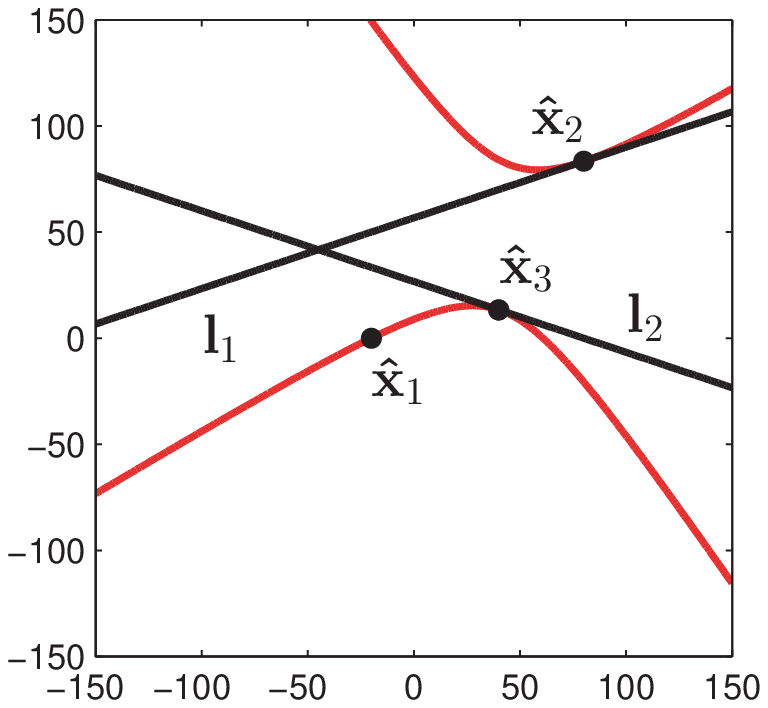}}
\subfigure[]{
\label{fig:subfig:b}  
\includegraphics[width=0.26\textwidth]{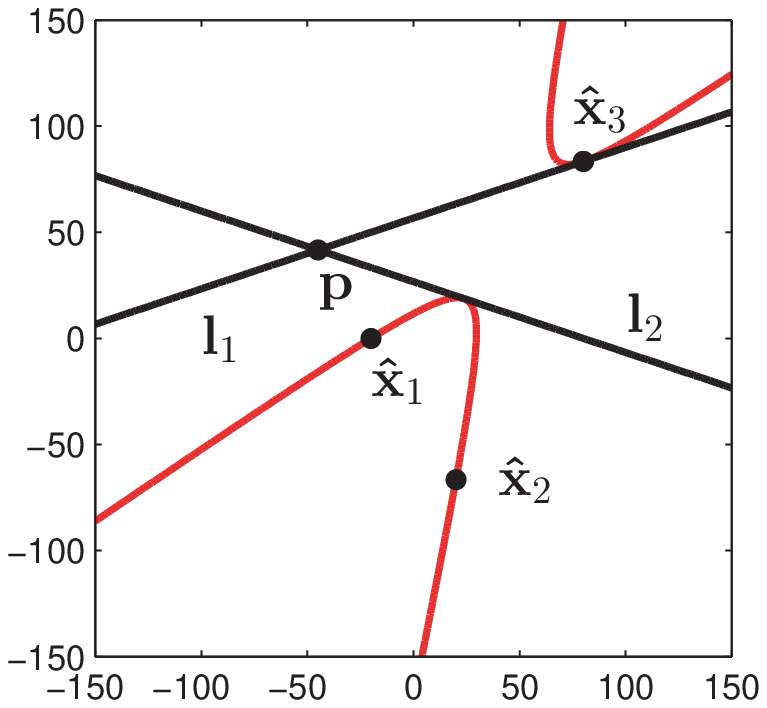}}
\subfigure[]{
\label{fig:subfig:c}  
\includegraphics[width=0.26\textwidth]{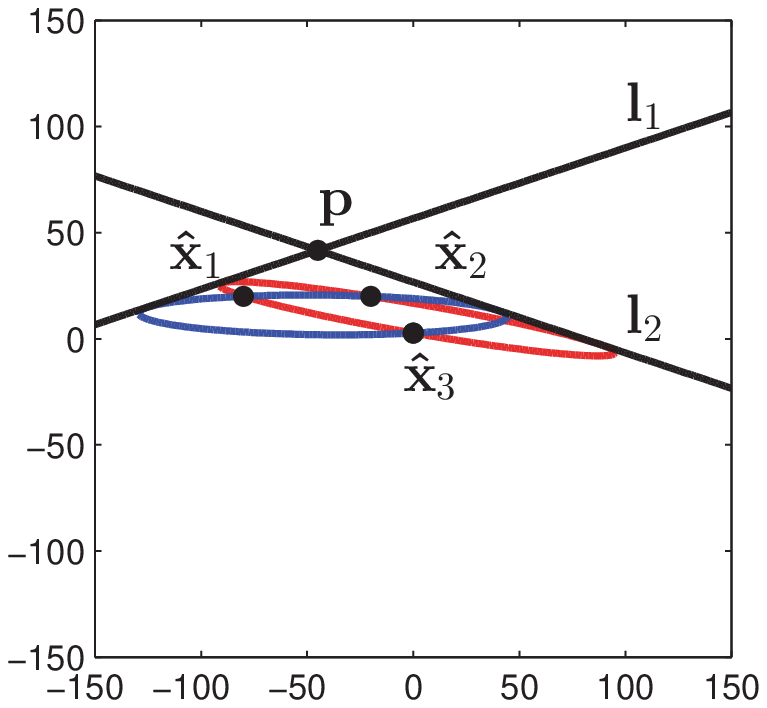}}
\vspace{1mm}
\subfigure[]{
\label{fig:subfig:a}  
\includegraphics[width=0.26\textwidth]{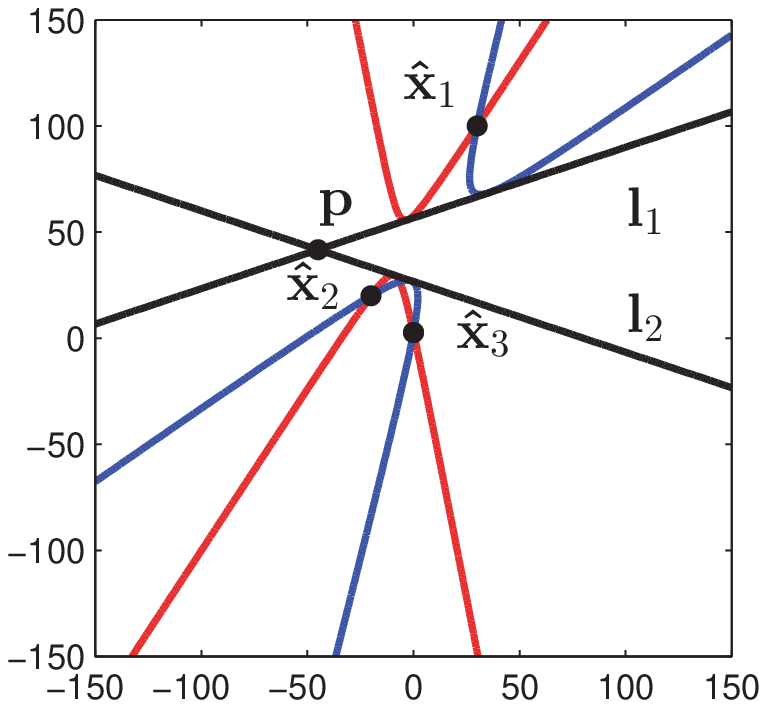}}
\subfigure[]{
\label{fig:subfig:b}  
\includegraphics[width=0.26\textwidth]{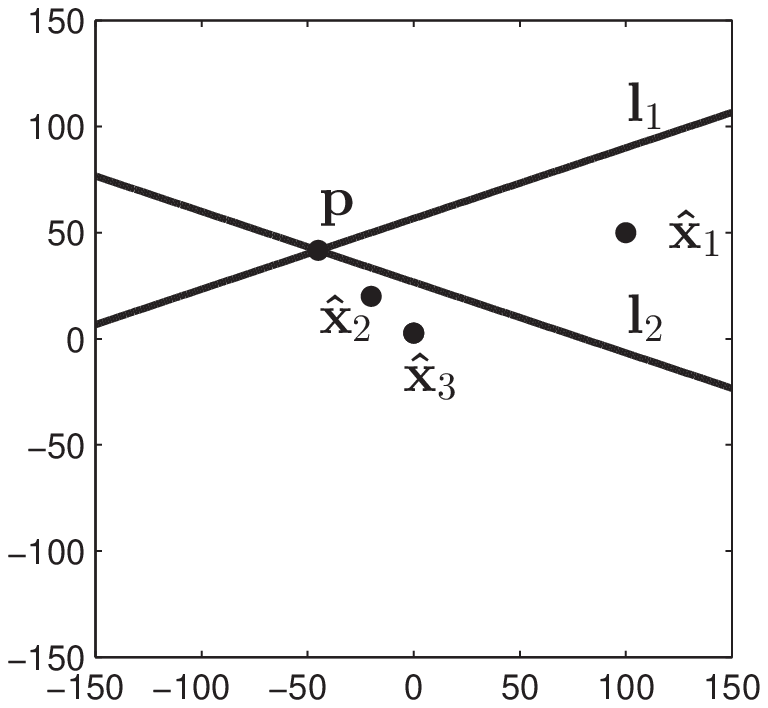}}
\subfigure[]{
\label{fig:subfig:c}  
\includegraphics[width=0.26\textwidth]{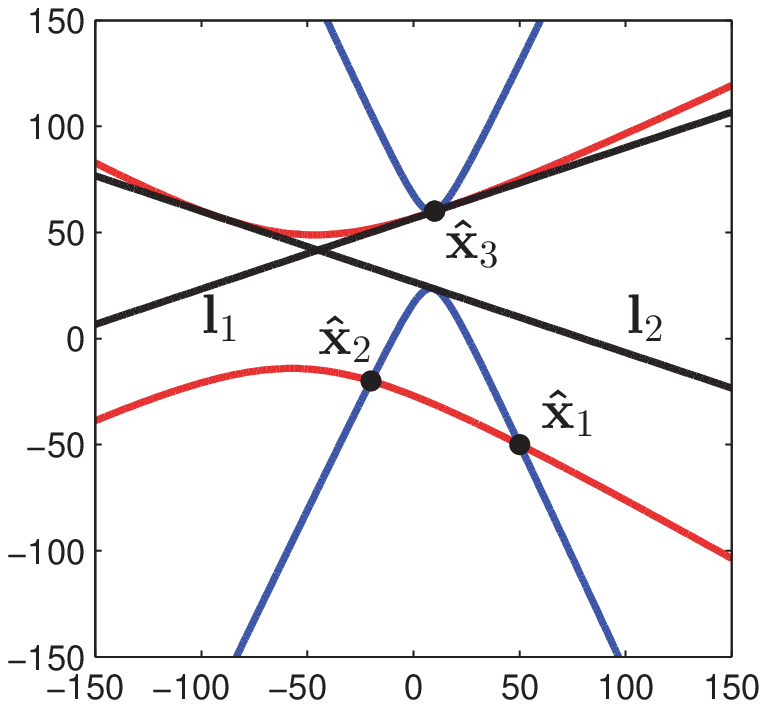}}
\vspace{1mm}
\subfigure[]{
\label{fig:subfig:a}  
\includegraphics[width=0.26\textwidth]{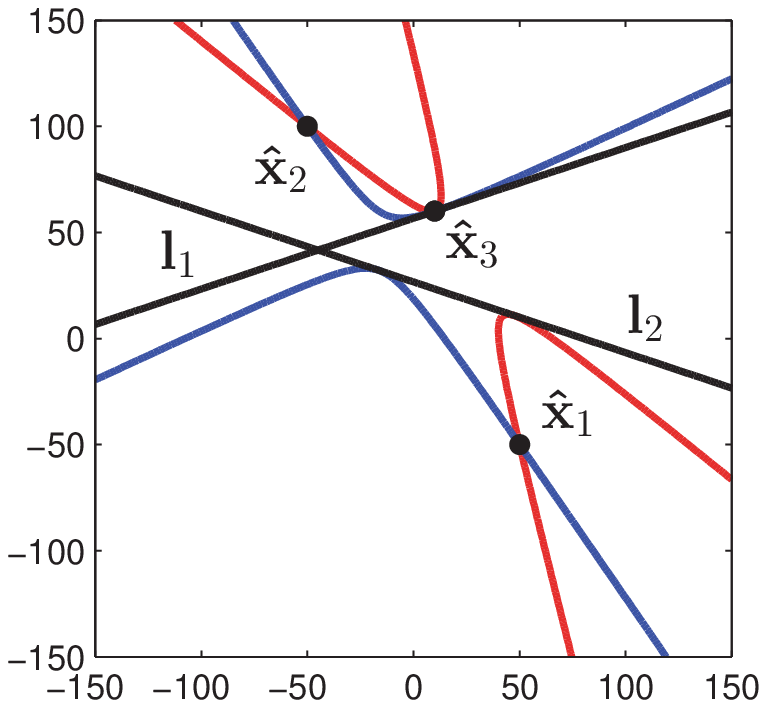}}
\subfigure[]{
\label{fig:subfig:b}  
\includegraphics[width=0.26\textwidth]{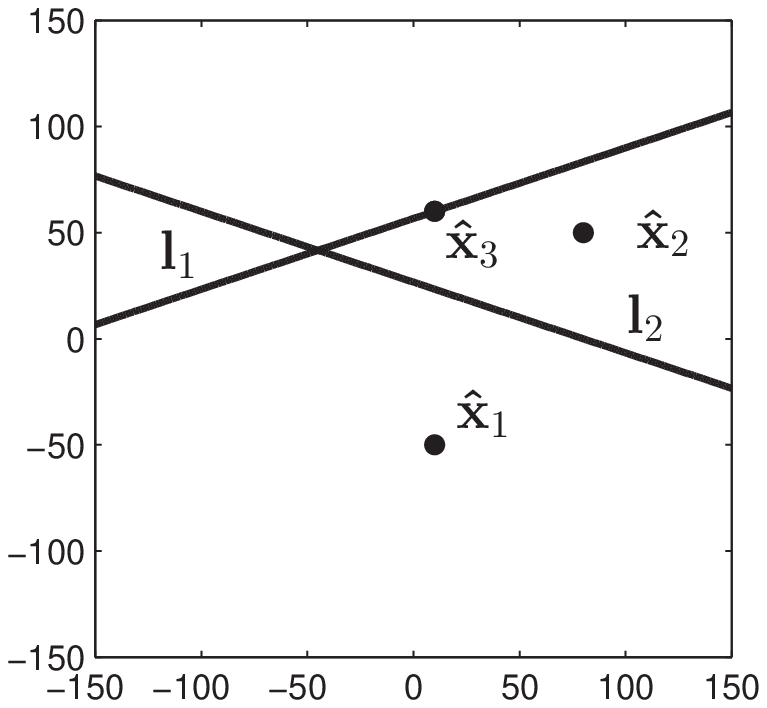}}
\subfigure[]{
\label{fig:subfig:c}  
\includegraphics[width=0.26\textwidth]{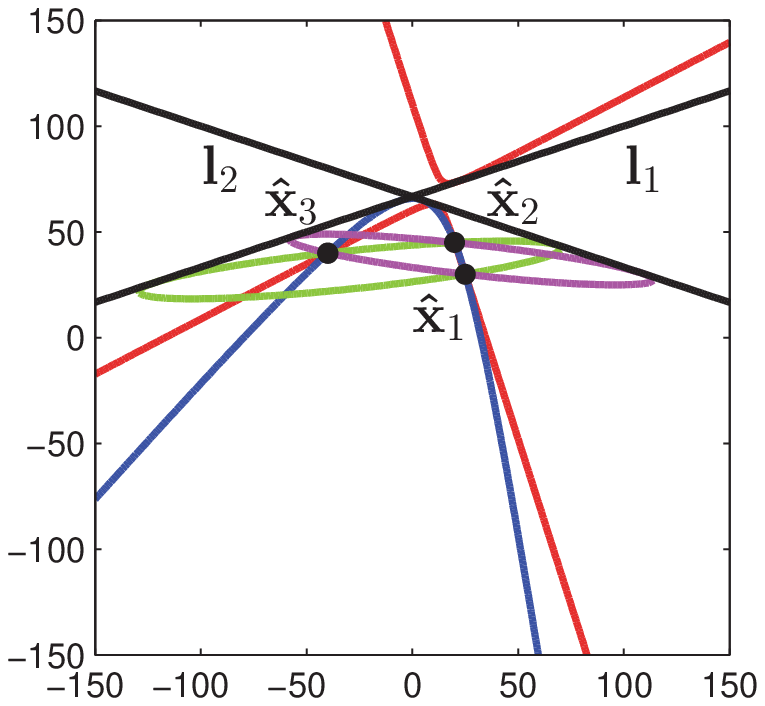}}
\vspace{1mm}
\subfigure[]{
\label{fig:subfig:a}  
\includegraphics[width=0.26\textwidth]{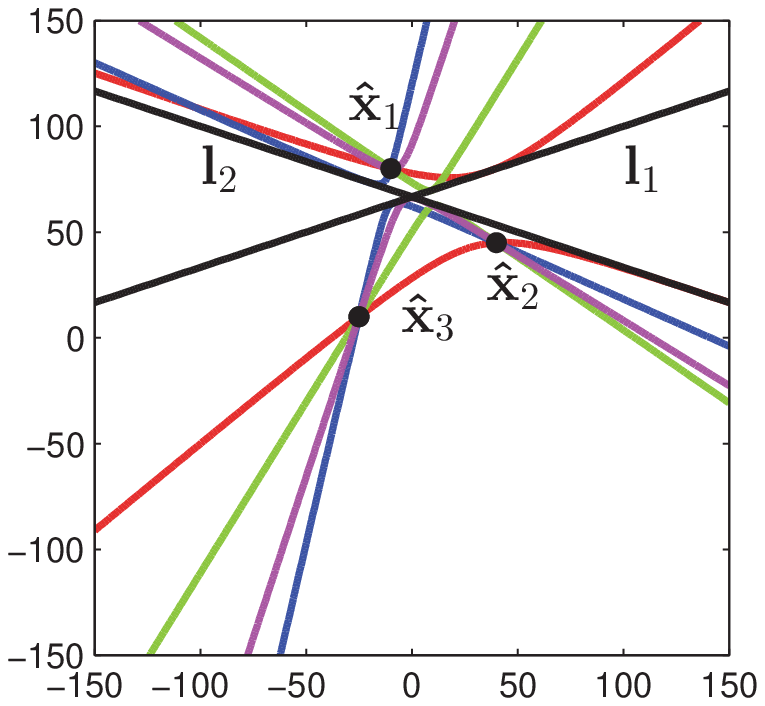}}
\subfigure[]{
\label{fig:subfig:b}  
\includegraphics[width=0.26\textwidth]{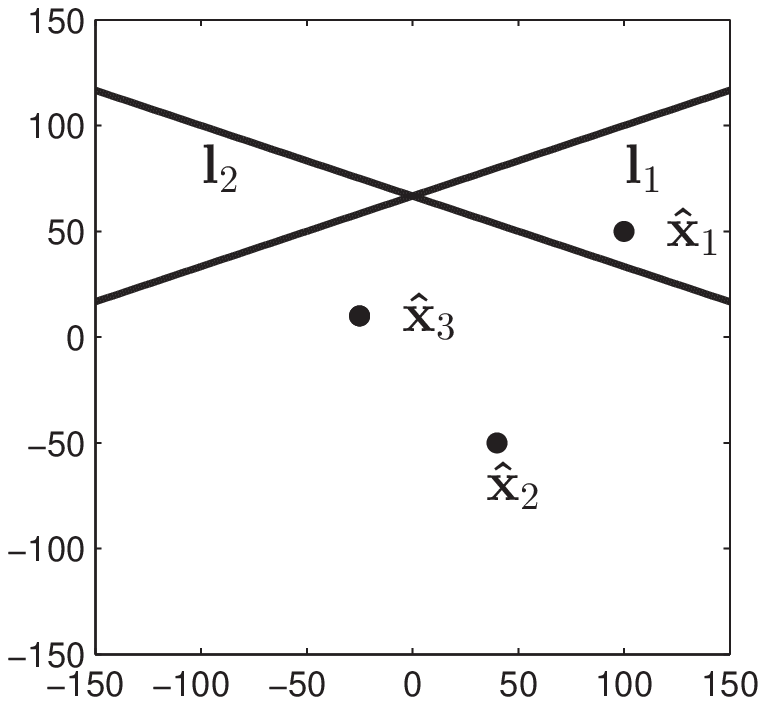}}
\subfigure[]{
\label{fig:subfig:c}  
\includegraphics[width=0.26\textwidth]{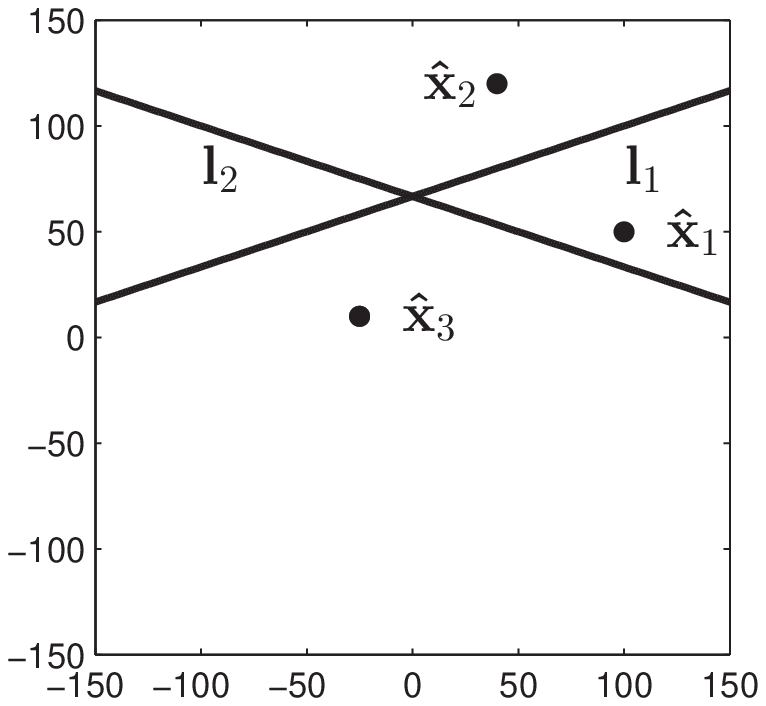}}
\caption{Conics from three points and two lines in general position. (a) shows a unique solution by Theorem 5 and Case 1 of Algorithm 3. (b) shows a unique solution by (ii) of Theorem 6 and Case 2 of Algorithm 3. (c) and (d) show two distinct real solutions by (i) of Theorem 6 and Case 3 of Algorithm 3. (e) shows no real solutions by Proposition 4. (f) and (g) show two distinct real solutions by Proposition 5 and Case 4 of Algorithm 3. (h) shows no real solutions by Proposition 5. (i) and (j) show four distinct real solutions by Proposition 3 and Case 5 of Algorithm 3. (k) and (l) show no real solutions by Proposition 3.}
\label{fig:subfig}  
\end{center}
\end{figure}

\section{Conclusion}

This paper presents a series of closed-form solutions to the conics from all minimal configurations of points and lines in general position. Using these closed-form formulae, we propose the algebraic and geometric conditions for a non-degenerate real conic from each of all minimal configurations. We also demonstrate the validity of our results by some examples. In fact, our closed-form solutions to the conics can be applied further into distinguishing type of a non-degenerate real conic from each of all minimal configurations. In future research, we will focus on the problems.\\

\noindent\textbf{Appendix A: If there exist two same eigenvalues, then $\textbf l_1\times\textbf l_2$ must lie on some side of the triangle $\hat {\textbf x}_1\hat {\textbf x}_2\hat {\textbf x}_3$.}\\

We prove the case of $\lambda_1=\lambda_2$, and the proofs of other cases are similar. If $\lambda_1=\lambda_2$, then we have $\hat {\textbf x}_2^T\textbf l_1=\rho\hat {\textbf x}_3^T\textbf l_1$ and $\hat {\textbf x}_2^T\textbf l_2=\rho\hat {\textbf x}_3^T\textbf l_2$, namely the point $\hat {\textbf x}_2-\rho\hat {\textbf x}_3$ does not only lie on the line $\textbf l_1$, but also the line $\textbf l_2$. That means that $\textbf l_1\times\textbf l_2\simeq\hat {\textbf x}_2-\rho\hat {\textbf x}_3$, namely $\textbf l_1\times\textbf l_2$ must lie on the side through $\hat {\textbf x}_2$, $\hat {\textbf x}_3$ of the triangle $\hat {\textbf x}_1\hat {\textbf x}_2\hat {\textbf x}_3$.\\




\begin{small}
\noindent\textbf{Acknowledgments}\quad The author thanks the anonymous referees for valuable suggestions.
\end{small}

\section*{References}
\begin{enumerate}
\item S. Kleiman, Chasles's enumerative theory of conics: A historical introduction, in studies in Algebraic Geometry, Mathematical Association of America Studies in Mathematics, Mathematical Association of America, Washington, DC (20)1980 117-138.
\item A. Bashelor, A. Ksir, W. Traves, Enumerative algebraic geometry of conics, American Mathematical Monthly 115(2008) 701-728.
\item H. D\"{o}rrie, 100 great problems of elementary mathematics: their history and solution - translated by David Antin, Dover Publications, New York, 1965.
\item E. Lord, Symmetry and Pattern in Projective Geometry, Springer-Verlag London, Springer, London, 2013.
\item R. I. Hartley, A. Zisserman, Multiple view geometry in computer vision, second ed., Cambridge, U.K.: Cambridge University Press, 2004.
\item F. Sottile, Real solutions to equations from geometry, University Lecture Series, American Mathematical Society, Providence, Rhode Island, 2011.
\item F. Sottile, The special Schubert calculus is real, Electron. Res. Announc. AMS 5(1999) 35¨C39.
\item M. A. Fishler, R. C. Bolles, Random sample consensus: a paradigm for model fitting with applications to image analysis and automated cartography. Comm. ACM. 24(1981) 381-395.
\item J. Semple, G. Kneebone, Algebraic Projective Geometry, Oxford University Press, 1979.
\item S. Axler, Linear Algebra Done Right, second ed., Springer-Verlag, New York, 1997.
\item M. Artin, Algebra, second ed., Pearson Education, Inc, Prentice Hall, 2010.
\item S. Roman, Advanced linear algebra, third ed., Springer, New York, 2008.
\end{enumerate}
\bibliographystyle{model2-names}
\bibliography{refs}
\end{document}